\documentclass[onefignum,onetabnum]{siamart171218}

\usepackage{lipsum}
\usepackage{amsfonts}
\usepackage{graphicx}
\usepackage{epstopdf}
\usepackage{algorithmic}
\usepackage{url,tikz,verbatim,booktabs,amsmath,bm,amssymb,subfigure}

\ifpdf
  \DeclareGraphicsExtensions{.eps,.pdf,.png,.jpg}
\else
  \DeclareGraphicsExtensions{.eps}
\fi

\newcommand{\mat}[1]{\bm{#1}}
\newcommand{\ten}[1]{\bm{\mathcal{#1}}}

\newsiamremark{remark}{Remark}
\newsiamremark{problem}{Problem}
\newsiamremark{hypothesis}{Hypothesis}
\crefname{hypothesis}{Hypothesis}{Hypotheses}
\newsiamthm{claim}{Claim}

\headers{MERACLE}{K. Batselier \textit{et al.}}

\title{MERACLE: Constructive layer-wise conversion of a Tensor Train into a MERA}

\author{Kim Batselier\thanks{Delft Center for Systems and Control, Delft University of Technology, Delft, The Netherlands.
  (\email{k.batselier@tudelft.nl}).} \and  Andrzej Cichocki\thanks{ Skolkovo Institute of Science and Technology
(Skoltech), Russia.} \and Ngai Wong\thanks{Department of Electrical and Electronic Engineering, The University of Hong Kong, Hong Kong.}}
  
\usepackage{amsopn}

\ifpdf
\hypersetup{
  pdftitle={TT2TM},
  pdfauthor={K. Batselier}
}
\fi

\begin{document}

\maketitle

% REQUIRED
\begin{abstract}
In this article two new algorithms are presented that convert a given data tensor train into either a Tucker decomposition with orthogonal matrix factors or a multi-scale entanglement renormalization ansatz (MERA). The Tucker core tensor is never explicitly computed but stored as a tensor train instead, resulting in both computationally and storage efficient algorithms. Both the multilinear Tucker-ranks as well as the MERA-ranks are automatically determined by the algorithm for a given upper bound on the relative approximation error. In addition, an iterative algorithm with low computational complexity based on solving an orthogonal Procrustes problem is proposed for the first time to retrieve optimal rank-lowering disentangler tensors, which are a crucial component in the construction of a low-rank MERA. Numerical experiments demonstrate the effectiveness of the proposed algorithms together with the potential storage benefit of a low-rank MERA over a tensor train.
\end{abstract}

% REQUIRED

% REQUIRED

\section{Introduction}
Tensor decompositions have played an important role over the past 2 decades in lifting the curse of dimensionality in myriad of applications~\cite{MAL-059,cichocki2015tensor,MAL-067,Khoromskij2011,sidiropoulos2017tensor}. The key idea in lifting the curse of dimensionality with tensor decompositions is the usage of a low-rank approximation. Many kinds of decompositions have consequently been developed and each has its own rank definition. The canonical polyadic decomposition (CPD)~\cite{candecomp,harshman1970fpp,Hitchcock1927} and Tucker decomposition~\cite{candecomp,tuckerreview} both generalize the notion of the matrix singular value decomposition (SVD) to higher order tensors and have therefore received a lot of attention. More recent tensor decompositions are the Tensor Train~\cite{ivanTT,Espig2011,Espig2012,Khoromskij2011} (TT) and hierarchical Tucker decomposition~\cite{Grasedyck2010,Hackbusch2009}. It turns out that the latter two decompositions were already known in the quantum mechanics and condensed matter physics communities as the matrix product state (MPS)~\cite{MPS1997} and Tensor Tree Network~\cite{TTN2006}, respectively. The multi-scale entanglement renormalization ansatz (MERA)~\cite{MERA2009,MERA2008} is an extension of the TTN decomposition, recently proposed in quantum mechanics but has so far not received enough attention in the numerical linear algebra community. A key component of the MERA is the so-called disentangler tensor, responsible for limiting the growth of the TTN-ranks over consecutive levels. Although the computation of a MERA from a given tensor can be deduced from~\cite{MERA2009}, computations are intensive due to multiple contractions and do not allow for the discovery of optimal ranks of the decomposition. The contributions of this article address this area. Specifically, we
\begin{enumerate}
    \item propose an algorithm that converts a given TT into a Tucker decomposition with guaranteed error bounds. 
    \item propose an algorithm that converts a given TT into a MERA with guaranteed error bounds. This algorithm is called MERA Constructive Layer-wise Expansion (MERACLE).
    \item propose an iterative algorithm that computes a rank-lowering disentangler.
\end{enumerate}
The resulting ranks of the computed Tucker and MERA approximations are completely determined by a given upper bound on the relative approximation error. The conversion of a TT into a Tucker decomposition was first suggested in~\cite{TTTucker}, where the corresponding algorithm uses an iterative Alternating Least Squares (ALS) approach. It will be shown in this article that no ALS procedure is necessary. In fact, for a $D$-th order tensor it is sufficient to perform $D$ consecutive SVD computations as described in Algorithm~\ref{alg:TT2Tucker}. It is then shown in Algorithm~\ref{alg:TT2MERA} that a TT can be converted into an $L$-layer MERA by applying Algorithm~\ref{alg:TT2Tucker} $2L$ times. The obtained MERA ranks are, however, not optimal and this is identified to be due to the disentangler tensor computation. An iterative orthogonal Procrustes algorithm is proposed that, to our knowledge for the first time ever, is able to compute optimal disentanglers that result in a minimal-rank MERA.

In Section~\ref{sec:basics} we introduce the notation and relevant tensor decompositions. The algorithm that converts a given TT into a Tucker decomposition with a guaranteed relative error bound is fully described in Section~\ref{sec:TT2Tucker}. The application of Algorithm~\ref{alg:TT2Tucker} for the conversion of a given TT into a MERA with a guaranteed relative error bound is illustrated in Section~\ref{sec:TT2MERA}. Section~\ref{sec:Procrustes} discusses the problem of finding optimal disentangler tensors and the iterative Procrustes algorithm is proposed. Finally, in Section~\ref{sec:experiments} numerical experiments demonstrate the effectiveness of the proposed algorithms.

\section{Tensor basics}
\label{sec:basics}
A $D$-way or $D$th order tensor $\ten{A} \in\mathbb{R}^{I_1\times I_2\times \cdots\times I_D}$ is a $D$-dimensional array where each entry is completely determined by $D$ indices $i_1,\ldots,i_D$. The scalar $D$ is also often called the order of the tensor. The convention $i_d=1,2,\ldots,I_d$ is used, together with MATLAB colon notation. Boldface capital calligraphic letters $\ten{A},\ten{B},\ldots$ are used to denote tensors, boldface capital letters $\mat{A},\mat{B},\ldots$ denote matrices, boldface letters $\mat{a},\mat{b},\ldots$ denote vectors, and Roman letters $a,b,\ldots$ denote scalars. The identity matrix of order $N$ is denoted $\mat{I}_N$. The Frobenius norm $||\ten{A}||_F^2$ of a tensor $\ten{A}$ is defined as the sum of squares of all tensor entries. The order of a tensor can be altered by grouping several indices together into a multi-index. The conversion of a multi-index $[i_1i_2\cdots i_D]$ into a linear index is per definition
\begin{align}
[i_1i_2\cdots i_D] := i_1 + \sum_{k=2}^{D}\,(i_k-1)\,\prod_{l=1}^{k-1} I_l.
\label{eq:multiindex}
\end{align}
In what follows, we will introduce three important tensor operations. The first tensor operation is the ``reshape" operation, which changes the order of a given tensor and is commonly used to flatten tensors into matrices and vice versa.
\begin{definition} The operator ``reshape($\ten{A},[J_1,J_2,\ldots,J_K])$" reshapes the $d$-way tensor $\ten{A} \in\mathbb{R}^{I_1\times I_2\times \cdots\times I_D}$ into a tensor with dimensions $J_1 \times J_2 \times \cdots \times J_K$, with $\prod_{d=1}^D I_d = \prod_{k=1}^K J_k$.
\end{definition}
Another important operation is the generalization of the matrix transpose to three or more indices.
\begin{definition}The operator ``permute($\ten{A},\mat{p})$" rearranges the indices of {$\ten{A} \in\mathbb{R}^{I_1\times I_2\times \cdots\times I_D}$} so that they are in the order specified by the vector $\mat{p}$. The resulting tensor has the same values of $\ten{A}$ but the order of the subscripts needed to access any particular element is rearranged as specified by $\mat{p}$. All the elements of $\mat{p}$ must be unique, real, positive, integer values from 1 to $D$.
\end{definition}
The definition of the ``permute" operation allows one to write the transpose of a matrix $\mat{A}$ as $\textrm{permute}(\mat{A},[2,1])$. By combining both the reshape and permute operations, we can now introduce the mode-$d$ matricization $\mat{A}_{<d>}$ of a tensor.
\begin{definition}(~\cite[p.~459]{tensorreview})
The mode-$d$ matricization $\mat{A}_{<d>}$ of a $D$-way tensor~$\ten{A}$ is the matrix with elements
\begin{align*}
\mat{A}_{<d>}( i_d, [i_1\cdots i_{d-1}i_{d+1}\cdots i_D]) &:= \ten{A}(i_1,i_2,\cdots,i_D).
\end{align*}
\end{definition}
The mode-$d$ matricization $\mat{A}_{<d>}$ is hence obtained from $\ten{A}$ as
\begin{align*}
   \mat{A}_{<d>} &= \textrm{reshape}(\textrm{permute}(\ten{A} ,[d,1,2,\ldots,d-1,d+1,\ldots,D]),[I_d,I_1\cdots I_D]). 
\end{align*}

The third and final important tensor operation is the summation over indices, also called contraction of indices. A particular common operation in this regard is the $d$-mode product of a tensor with a matrix.
\begin{definition}(~\cite[p.~460]{tensorreview})
The $d$-mode product, denoted $\ten{A} \times_d \mat{U}_d$, of a tensor $\ten{A} \in \mathbb{R}^{I_1 \times \cdots \times I_D}$ with a matrix $\mat{U}_d \in \mathbb{R}^{S_d \times I_d}$ is the tensor $\ten{B} \in \mathbb{R}^{I_1 \times I_{d-1} \times S_d \times I_{d+1} \times \cdots \times I_D}$ with elements
\begin{align*}
    \ten{B}(i_1,\ldots,i_{d-1},s_d,i_{d+1},\ldots,i_D) &:= \sum_{i_d=1}^{I_d} \ten{A}(i_1,\ldots,i_{d-1},i_d,i_{d+1},\ldots,i_D) \; \mat{U}_{d}(j_d,i_d).
\end{align*}
\end{definition}

\begin{figure}[t]
\subfigure[Diagram representation of a scalar $a$, vector $\mat{a}$, matrix $\mat{A}$ and 3-way tensor $\ten{A}$.]{
\begin{minipage}[b]{0.45\linewidth}
\centering
% Graphic for TeX using PGF
% Title: /home/kim/Dropbox/Work/Papers/Journal/16.MERA/figs/TNgraphs.dia
% Creator: Dia v0.97.2
% CreationDate: Wed Apr 19 13:58:56 2017
% For: kim
% \usepackage{tikz}
% The following commands are not supported in PSTricks at present
% We define them conditionally, so when they are implemented,
% this pgf file will use them.
\ifx\du\undefined
  \newlength{\du}
\fi
\setlength{\du}{4.5\unitlength}
\begin{tikzpicture}
\pgftransformxscale{1.000000}
\pgftransformyscale{-1.000000}
\definecolor{dialinecolor}{rgb}{0.000000, 0.000000, 0.000000}
\pgfsetstrokecolor{dialinecolor}
\definecolor{dialinecolor}{rgb}{1.000000, 1.000000, 1.000000}
\pgfsetfillcolor{dialinecolor}
\definecolor{dialinecolor}{rgb}{1.000000, 1.000000, 1.000000}
\pgfsetfillcolor{dialinecolor}
\pgfpathellipse{\pgfpoint{-6.073446\du}{10.779091\du}}{\pgfpoint{2.900000\du}{0\du}}{\pgfpoint{0\du}{2.800000\du}}f
\pgfusepath{fill}
\pgfsetlinewidth{0.100000\du}
\pgfsetdash{}{0pt}
\pgfsetdash{}{0pt}
\definecolor{dialinecolor}{rgb}{0.000000, 0.000000, 0.000000}
\pgfsetstrokecolor{dialinecolor}
\pgfpathellipse{\pgfpoint{-6.073446\du}{10.779091\du}}{\pgfpoint{2.900000\du}{0\du}}{\pgfpoint{0\du}{2.800000\du}}
\pgfusepath{stroke}
\pgfsetlinewidth{0.100000\du}
\pgfsetdash{}{0pt}
\pgfsetdash{}{0pt}
\pgfsetbuttcap
{
\definecolor{dialinecolor}{rgb}{0.000000, 0.000000, 0.000000}
\pgfsetfillcolor{dialinecolor}
% was here!!!
\definecolor{dialinecolor}{rgb}{0.000000, 0.000000, 0.000000}
\pgfsetstrokecolor{dialinecolor}
\draw (-4.022836\du,12.758990\du)--(-4.022836\du,16.80\du);
}
\pgfsetlinewidth{0.100000\du}
\pgfsetdash{}{0pt}
\pgfsetdash{}{0pt}
\pgfsetbuttcap
{
\definecolor{dialinecolor}{rgb}{0.000000, 0.000000, 0.000000}
\pgfsetfillcolor{dialinecolor}
% was here!!!
\definecolor{dialinecolor}{rgb}{0.000000, 0.000000, 0.000000}
\pgfsetstrokecolor{dialinecolor}
\draw (-8.124055\du,12.758990\du)--(-8.124055\du,16.80\du);
}
\pgfsetlinewidth{0.100000\du}
\pgfsetdash{}{0pt}
\pgfsetdash{}{0pt}
\pgfsetbuttcap
{
\definecolor{dialinecolor}{rgb}{0.000000, 0.000000, 0.000000}
\pgfsetfillcolor{dialinecolor}
% was here!!!
\definecolor{dialinecolor}{rgb}{0.000000, 0.000000, 0.000000}
\pgfsetstrokecolor{dialinecolor}
\draw (-6.073446\du,13.579091\du)--(-6.073446\du,16.80\du);
}
% setfont left to latex
\definecolor{dialinecolor}{rgb}{0.000000, 0.000000, 0.000000}
\pgfsetstrokecolor{dialinecolor}
\node[anchor=west] at (-8.20\du,10.6\du){$\ten{A}$};
\definecolor{dialinecolor}{rgb}{1.000000, 1.000000, 1.000000}
\pgfsetfillcolor{dialinecolor}
\pgfpathellipse{\pgfpoint{-36.297515\du}{10.779091\du}}{\pgfpoint{2.900000\du}{0\du}}{\pgfpoint{0\du}{2.800000\du}}
\pgfusepath{fill}
\pgfsetlinewidth{0.100000\du}
\pgfsetdash{}{0pt}
\pgfsetdash{}{0pt}
\definecolor{dialinecolor}{rgb}{0.000000, 0.000000, 0.000000}
\pgfsetstrokecolor{dialinecolor}
\pgfpathellipse{\pgfpoint{-36.297515\du}{10.779091\du}}{\pgfpoint{2.900000\du}{0\du}}{\pgfpoint{0\du}{2.800000\du}}
\pgfusepath{stroke}
% setfont left to latex
\definecolor{dialinecolor}{rgb}{0.000000, 0.000000, 0.000000}
\pgfsetstrokecolor{dialinecolor}
\node[anchor=west] at (-37.80\du,10.80\du){$a$};
\definecolor{dialinecolor}{rgb}{1.000000, 1.000000, 1.000000}
\pgfsetfillcolor{dialinecolor}
\pgfpathellipse{\pgfpoint{-26.560105\du}{10.779091\du}}{\pgfpoint{2.900000\du}{0\du}}{\pgfpoint{0\du}{2.800000\du}}
\pgfusepath{fill}
\pgfsetlinewidth{0.100000\du}
\pgfsetdash{}{0pt}
\pgfsetdash{}{0pt}
\definecolor{dialinecolor}{rgb}{0.000000, 0.000000, 0.000000}
\pgfsetstrokecolor{dialinecolor}
\pgfpathellipse{\pgfpoint{-26.560105\du}{10.779091\du}}{\pgfpoint{2.900000\du}{0\du}}{\pgfpoint{0\du}{2.800000\du}}
\pgfusepath{stroke}
\pgfsetlinewidth{0.100000\du}
\pgfsetdash{}{0pt}
\pgfsetdash{}{0pt}
\pgfsetbuttcap
{
\definecolor{dialinecolor}{rgb}{0.000000, 0.000000, 0.000000}
\pgfsetfillcolor{dialinecolor}
% was here!!!
\definecolor{dialinecolor}{rgb}{0.000000, 0.000000, 0.000000}
\pgfsetstrokecolor{dialinecolor}
\draw (-26.560105\du,13.579091\du)--(-26.560105\du,16.80\du);
}
% setfont left to latex
\definecolor{dialinecolor}{rgb}{0.000000, 0.000000, 0.000000}
\pgfsetstrokecolor{dialinecolor}
\node[anchor=west] at (-28.20\du,10.80\du){$\mat{a}$};
\definecolor{dialinecolor}{rgb}{1.000000, 1.000000, 1.000000}
\pgfsetfillcolor{dialinecolor}
\pgfpathellipse{\pgfpoint{-16.948419\du}{10.779091\du}}{\pgfpoint{2.900000\du}{0\du}}{\pgfpoint{0\du}{2.800000\du}}
\pgfusepath{fill}
\pgfsetlinewidth{0.100000\du}
\pgfsetdash{}{0pt}
\pgfsetdash{}{0pt}
\definecolor{dialinecolor}{rgb}{0.000000, 0.000000, 0.000000}
\pgfsetstrokecolor{dialinecolor}
\pgfpathellipse{\pgfpoint{-16.948419\du}{10.779091\du}}{\pgfpoint{2.900000\du}{0\du}}{\pgfpoint{0\du}{2.800000\du}}
\pgfusepath{stroke}
\pgfsetlinewidth{0.100000\du}
\pgfsetdash{}{0pt}
\pgfsetdash{}{0pt}
\pgfsetbuttcap
{
\definecolor{dialinecolor}{rgb}{0.000000, 0.000000, 0.000000}
\pgfsetfillcolor{dialinecolor}
% was here!!!
\definecolor{dialinecolor}{rgb}{0.000000, 0.000000, 0.000000}
\pgfsetstrokecolor{dialinecolor}
\draw (-14.897810\du,12.758990\du)--(-14.897810\du,16.80\du);
}
\pgfsetlinewidth{0.100000\du}
\pgfsetdash{}{0pt}
\pgfsetdash{}{0pt}
\pgfsetbuttcap
{
\definecolor{dialinecolor}{rgb}{0.000000, 0.000000, 0.000000}
\pgfsetfillcolor{dialinecolor}
% was here!!!
\definecolor{dialinecolor}{rgb}{0.000000, 0.000000, 0.000000}
\pgfsetstrokecolor{dialinecolor}
\draw (-18.999029\du,12.758990\du)--(-18.999029\du,16.80\du);
}
% setfont left to latex
\definecolor{dialinecolor}{rgb}{0.000000, 0.000000, 0.000000}
\pgfsetstrokecolor{dialinecolor}
\node[anchor=west] at (-19.00\du,10.6\du){$\mat{A}$};
\end{tikzpicture}
\end{minipage}}
\subfigure[Diagram representation of equation~\eqref{eqn:simpleTN} with all dimensions labelled.]{
\begin{minipage}[b]{0.45\linewidth}
\centering
  % Graphic for TeX using PGF
% Title: /home/kbatselier/Pictures/Dia/TNcontraction.dia
% Creator: Dia v0.97+git
% CreationDate: Wed Mar 13 15:59:08 2019
% For: kbatselier
% \usepackage{tikz}
% The following commands are not supported in PSTricks at present
% We define them conditionally, so when they are implemented,
% this pgf file will use them.
\ifx\du\undefined
  \newlength{\du}
\fi
\setlength{\du}{4\unitlength}
\begin{tikzpicture}[even odd rule]
\pgftransformxscale{1.000000}
\pgftransformyscale{-1.000000}
\definecolor{dialinecolor}{rgb}{0.000000, 0.000000, 0.000000}
\pgfsetstrokecolor{dialinecolor}
\pgfsetstrokeopacity{1.000000}
\definecolor{diafillcolor}{rgb}{1.000000, 1.000000, 1.000000}
\pgfsetfillcolor{diafillcolor}
\pgfsetfillopacity{1.000000}
\pgfsetlinewidth{0.200000\du}
\pgfsetdash{}{0pt}
\pgfsetbuttcap
{
\definecolor{diafillcolor}{rgb}{0.000000, 0.000000, 0.000000}
\pgfsetfillcolor{diafillcolor}
\pgfsetfillopacity{1.000000}
% was here!!!
\definecolor{dialinecolor}{rgb}{0.000000, 0.000000, 0.000000}
\pgfsetstrokecolor{dialinecolor}
\pgfsetstrokeopacity{1.000000}
\draw (17.845000\du,8.660000\du)--(13.345000\du,8.610000\du);
}
\pgfsetlinewidth{0.200000\du}
\pgfsetdash{}{0pt}
\pgfsetbuttcap
{
\definecolor{diafillcolor}{rgb}{0.000000, 0.000000, 0.000000}
\pgfsetfillcolor{diafillcolor}
\pgfsetfillopacity{1.000000}
% was here!!!
\definecolor{dialinecolor}{rgb}{0.000000, 0.000000, 0.000000}
\pgfsetstrokecolor{dialinecolor}
\pgfsetstrokeopacity{1.000000}
\draw (20.745000\du,11.460000\du)--(20.750000\du,16.929200\du);
}
\pgfsetlinewidth{0.100000\du}
\pgfsetdash{}{0pt}
\definecolor{diafillcolor}{rgb}{1.000000, 1.000000, 1.000000}
\pgfsetfillcolor{diafillcolor}
\pgfsetfillopacity{1.000000}
\pgfpathellipse{\pgfpoint{20.745000\du}{8.610000\du}}{\pgfpoint{2.900000\du}{0\du}}{\pgfpoint{0\du}{2.800000\du}}
\pgfusepath{fill}
\definecolor{dialinecolor}{rgb}{0.000000, 0.000000, 0.000000}
\pgfsetstrokecolor{dialinecolor}
\pgfsetstrokeopacity{1.000000}
\pgfpathellipse{\pgfpoint{20.745000\du}{8.610000\du}}{\pgfpoint{2.900000\du}{0\du}}{\pgfpoint{0\du}{2.800000\du}}
\pgfusepath{stroke}
% setfont left to latex
\definecolor{dialinecolor}{rgb}{0.000000, 0.000000, 0.000000}
\pgfsetstrokecolor{dialinecolor}
\pgfsetstrokeopacity{1.000000}
\definecolor{diafillcolor}{rgb}{0.000000, 0.000000, 0.000000}
\pgfsetfillcolor{diafillcolor}
\pgfsetfillopacity{1.000000}
\node[anchor=base west,inner sep=0pt,outer sep=0pt,color=dialinecolor] at (19.807500\du,9.604930\du){$\ten{A}$};
\pgfsetlinewidth{0.100000\du}
\pgfsetdash{}{0pt}
\definecolor{diafillcolor}{rgb}{1.000000, 1.000000, 1.000000}
\pgfsetfillcolor{diafillcolor}
\pgfsetfillopacity{1.000000}
\pgfpathellipse{\pgfpoint{10.445000\du}{8.610000\du}}{\pgfpoint{2.900000\du}{0\du}}{\pgfpoint{0\du}{2.800000\du}}
\pgfusepath{fill}
\definecolor{dialinecolor}{rgb}{0.000000, 0.000000, 0.000000}
\pgfsetstrokecolor{dialinecolor}
\pgfsetstrokeopacity{1.000000}
\pgfpathellipse{\pgfpoint{10.445000\du}{8.610000\du}}{\pgfpoint{2.900000\du}{0\du}}{\pgfpoint{0\du}{2.800000\du}}
\pgfusepath{stroke}
% setfont left to latex
\definecolor{dialinecolor}{rgb}{0.000000, 0.000000, 0.000000}
\pgfsetstrokecolor{dialinecolor}
\pgfsetstrokeopacity{1.000000}
\definecolor{diafillcolor}{rgb}{0.000000, 0.000000, 0.000000}
\pgfsetfillcolor{diafillcolor}
\pgfsetfillopacity{1.000000}
\node[anchor=base west,inner sep=0pt,outer sep=0pt,color=dialinecolor] at (9.507490\du,9.504930\du){$\mat{U}_1$};
\pgfsetlinewidth{0.200000\du}
\pgfsetdash{}{0pt}
\pgfsetbuttcap
{
\definecolor{diafillcolor}{rgb}{0.000000, 0.000000, 0.000000}
\pgfsetfillcolor{diafillcolor}
\pgfsetfillopacity{1.000000}
% was here!!!
\definecolor{dialinecolor}{rgb}{0.000000, 0.000000, 0.000000}
\pgfsetstrokecolor{dialinecolor}
\pgfsetstrokeopacity{1.000000}
\draw (7.545000\du,8.610000\du)--(2.945550\du,8.660550\du);
}
\pgfsetlinewidth{0.200000\du}
\pgfsetdash{}{0pt}
\pgfsetbuttcap
{
\definecolor{diafillcolor}{rgb}{0.000000, 0.000000, 0.000000}
\pgfsetfillcolor{diafillcolor}
\pgfsetfillopacity{1.000000}
% was here!!!
\definecolor{dialinecolor}{rgb}{0.000000, 0.000000, 0.000000}
\pgfsetstrokecolor{dialinecolor}
\pgfsetstrokeopacity{1.000000}
\draw (28.645000\du,8.610000\du)--(23.645000\du,8.660000\du);
}
\pgfsetlinewidth{0.100000\du}
\pgfsetdash{}{0pt}
\definecolor{diafillcolor}{rgb}{1.000000, 1.000000, 1.000000}
\pgfsetfillcolor{diafillcolor}
\pgfsetfillopacity{1.000000}
\pgfpathellipse{\pgfpoint{31.545000\du}{8.610000\du}}{\pgfpoint{2.900000\du}{0\du}}{\pgfpoint{0\du}{2.800000\du}}
\pgfusepath{fill}
\definecolor{dialinecolor}{rgb}{0.000000, 0.000000, 0.000000}
\pgfsetstrokecolor{dialinecolor}
\pgfsetstrokeopacity{1.000000}
\pgfpathellipse{\pgfpoint{31.545000\du}{8.610000\du}}{\pgfpoint{2.900000\du}{0\du}}{\pgfpoint{0\du}{2.800000\du}}
\pgfusepath{stroke}
% setfont left to latex
\definecolor{dialinecolor}{rgb}{0.000000, 0.000000, 0.000000}
\pgfsetstrokecolor{dialinecolor}
\pgfsetstrokeopacity{1.000000}
\definecolor{diafillcolor}{rgb}{0.000000, 0.000000, 0.000000}
\pgfsetfillcolor{diafillcolor}
\pgfsetfillopacity{1.000000}
\node[anchor=base west,inner sep=0pt,outer sep=0pt,color=dialinecolor] at (30.857500\du,9.304930\du){$\mat{u}_3$};
% setfont left to latex
\definecolor{dialinecolor}{rgb}{0.000000, 0.000000, 0.000000}
\pgfsetstrokecolor{dialinecolor}
\pgfsetstrokeopacity{1.000000}
\definecolor{diafillcolor}{rgb}{0.000000, 0.000000, 0.000000}
\pgfsetfillcolor{diafillcolor}
\pgfsetfillopacity{1.000000}
\node[anchor=base west,inner sep=0pt,outer sep=0pt,color=dialinecolor] at (3.925000\du,7.375000\du){$J_1$};
% setfont left to latex
\definecolor{dialinecolor}{rgb}{0.000000, 0.000000, 0.000000}
\pgfsetstrokecolor{dialinecolor}
\pgfsetstrokeopacity{1.000000}
\definecolor{diafillcolor}{rgb}{0.000000, 0.000000, 0.000000}
\pgfsetfillcolor{diafillcolor}
\pgfsetfillopacity{1.000000}
\node[anchor=base west,inner sep=0pt,outer sep=0pt,color=dialinecolor] at (14.750000\du,7.292812\du){$I_1$};
% setfont left to latex
\definecolor{dialinecolor}{rgb}{0.000000, 0.000000, 0.000000}
\pgfsetstrokecolor{dialinecolor}
\pgfsetstrokeopacity{1.000000}
\definecolor{diafillcolor}{rgb}{0.000000, 0.000000, 0.000000}
\pgfsetfillcolor{diafillcolor}
\pgfsetfillopacity{1.000000}
\node[anchor=base west,inner sep=0pt,outer sep=0pt,color=dialinecolor] at (21.425000\du,15.435000\du){$I_2$};
% setfont left to latex
\definecolor{dialinecolor}{rgb}{0.000000, 0.000000, 0.000000}
\pgfsetstrokecolor{dialinecolor}
\pgfsetstrokeopacity{1.000000}
\definecolor{diafillcolor}{rgb}{0.000000, 0.000000, 0.000000}
\pgfsetfillcolor{diafillcolor}
\pgfsetfillopacity{1.000000}
\node[anchor=base west,inner sep=0pt,outer sep=0pt,color=dialinecolor] at (25.200000\du,7.292812\du){$I_3$};
\end{tikzpicture}
\end{minipage}}
\caption{Basic TN diagrams.}
\label{fig:TNdiagrams}
\end{figure}
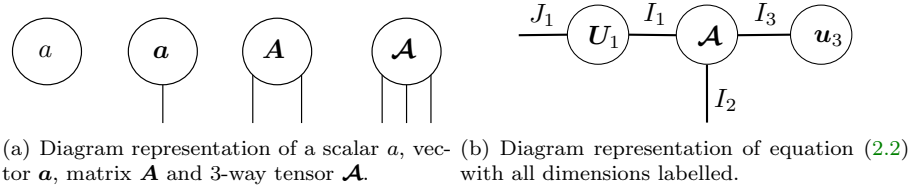

A very convenient graphical representation of $D$-way tensors is shown in Figure~\ref{fig:TNdiagrams}(a). Tensors are here represented by nodes and each edge denotes a particular index of the tensor. The order of the tensor is then easily determined by counting the number of edges. Since a scalar is a zeroth-order tensor, it is represented by a node without any edges.
The graphical representation of a summation over an index is by connecting the edge between the two nodes in the diagram. For example, the two index summations of a 3-way tensor $\ten{A} \in \mathbb{R}^{I_1 \times I_2 \times I_3}$ with a matrix $\mat{U}_1 \in \mathbb{R}^{J_1 \times I_1}$ and a vector $\mat{u}_3 \in \mathbb{R}^{I_3}$
\begin{align}
\label{eqn:simpleTN}
(\ten{A} \, \times_1\, \mat{U}_1 \, \times_3 \,\mat{u}_3^T)= \sum_{i_1,i_3} \, \ten{A}(i_1,:,i_3) \; \mat{U}_1(:,i_1) \; \mat{u}_3(i_3)
\end{align}
is graphically depicted in Figure~\ref{fig:TNdiagrams}(b) by two connected edges between the nodes for $\ten{A},\mat{U}_1$ and $\mat{u}_3$. The result from these two summations is a $J_1 \times I_2$ matrix, which can also be deduced from the two ``free" edges in Figure~\ref{fig:TNdiagrams}(b). Three important tensor decompositions in this article are the Tucker decomposition, the TT and the MERA. Each of these decompositions will now be briefly discussed.
\subsection{Tucker decomposition}
\begin{figure}[tb]
\subfigure[Diagram of a Tucker decomposition.]{
\begin{minipage}[b]{0.4\linewidth}
\centering
% Graphic for TeX using PGF
% Title: /home/kbatselier/Pictures/Dia/Tucker.dia
% Creator: Dia v0.97.3
% CreationDate: Tue Nov 20 11:54:37 2018
% For: kbatselier
% \usepackage{tikz}
% The following commands are not supported in PSTricks at present
% We define them conditionally, so when they are implemented,
% this pgf file will use them.
\ifx\du\undefined
  \newlength{\du}
\fi
\setlength{\du}{3\unitlength}
\begin{tikzpicture}
\pgftransformxscale{1.000000}
\pgftransformyscale{-1.000000}
\definecolor{dialinecolor}{rgb}{0.000000, 0.000000, 0.000000}
\pgfsetstrokecolor{dialinecolor}
\definecolor{dialinecolor}{rgb}{1.000000, 1.000000, 1.000000}
\pgfsetfillcolor{dialinecolor}
\definecolor{dialinecolor}{rgb}{1.000000, 1.000000, 1.000000}
\pgfsetfillcolor{dialinecolor}
\pgfpathellipse{\pgfpoint{31.725000\du}{9.575000\du}}{\pgfpoint{16.775000\du}{0\du}}{\pgfpoint{0\du}{5.075000\du}}
\pgfusepath{fill}
\pgfsetlinewidth{0.200000\du}
\pgfsetdash{}{0pt}
\pgfsetdash{}{0pt}
\definecolor{dialinecolor}{rgb}{0.000000, 0.000000, 0.000000}
\pgfsetstrokecolor{dialinecolor}
\pgfpathellipse{\pgfpoint{31.725000\du}{9.575000\du}}{\pgfpoint{16.775000\du}{0\du}}{\pgfpoint{0\du}{5.075000\du}}
\pgfusepath{stroke}
\pgfsetlinewidth{0.200000\du}
\pgfsetdash{}{0pt}
\pgfsetdash{}{0pt}
\pgfsetbuttcap
{
\definecolor{dialinecolor}{rgb}{0.000000, 0.000000, 0.000000}
\pgfsetfillcolor{dialinecolor}
% was here!!!
\definecolor{dialinecolor}{rgb}{0.000000, 0.000000, 0.000000}
\pgfsetstrokecolor{dialinecolor}
\draw (19.500000\du,13.150000\du)--(19.545554\du,20.005605\du);
}
\pgfsetlinewidth{0.200000\du}
\pgfsetdash{}{0pt}
\pgfsetdash{}{0pt}
\pgfsetbuttcap
{
\definecolor{dialinecolor}{rgb}{0.000000, 0.000000, 0.000000}
\pgfsetfillcolor{dialinecolor}
% was here!!!
\definecolor{dialinecolor}{rgb}{0.000000, 0.000000, 0.000000}
\pgfsetstrokecolor{dialinecolor}
\draw (27.717003\du,14.503783\du)--(27.775005\du,20.005973\du);
}
\definecolor{dialinecolor}{rgb}{1.000000, 1.000000, 1.000000}
\pgfsetfillcolor{dialinecolor}
\pgfpathellipse{\pgfpoint{19.571634\du}{23.930519\du}}{\pgfpoint{3.725000\du}{0\du}}{\pgfpoint{0\du}{3.825000\du}}
\pgfusepath{fill}
\pgfsetlinewidth{0.200000\du}
\pgfsetdash{}{0pt}
\pgfsetdash{}{0pt}
\definecolor{dialinecolor}{rgb}{0.000000, 0.000000, 0.000000}
\pgfsetstrokecolor{dialinecolor}
\pgfpathellipse{\pgfpoint{19.571634\du}{23.930519\du}}{\pgfpoint{3.725000\du}{0\du}}{\pgfpoint{0\du}{3.825000\du}}
\pgfusepath{stroke}
% setfont left to latex
\definecolor{dialinecolor}{rgb}{0.000000, 0.000000, 0.000000}
\pgfsetstrokecolor{dialinecolor}
\node[anchor=west] at (30.881242\du,10.151438\du){$\ten{S}$};
% setfont left to latex
\definecolor{dialinecolor}{rgb}{0.000000, 0.000000, 0.000000}
\pgfsetstrokecolor{dialinecolor}
\node[anchor=west] at (16\du,24.499956\du){$\mat{U}_1$};
\pgfsetlinewidth{0.200000\du}
\pgfsetdash{}{0pt}
\pgfsetdash{}{0pt}
\pgfsetbuttcap
{
\definecolor{dialinecolor}{rgb}{0.000000, 0.000000, 0.000000}
\pgfsetfillcolor{dialinecolor}
% was here!!!
\definecolor{dialinecolor}{rgb}{0.000000, 0.000000, 0.000000}
\pgfsetstrokecolor{dialinecolor}
\draw (19.571634\du,27.755519\du)--(19.505497\du,32.442013\du);
}
\definecolor{dialinecolor}{rgb}{1.000000, 1.000000, 1.000000}
\pgfsetfillcolor{dialinecolor}
\pgfpathellipse{\pgfpoint{27.816377\du}{23.930519\du}}{\pgfpoint{3.725000\du}{0\du}}{\pgfpoint{0\du}{3.825000\du}}
\pgfusepath{fill}
\pgfsetlinewidth{0.200000\du}
\pgfsetdash{}{0pt}
\pgfsetdash{}{0pt}
\definecolor{dialinecolor}{rgb}{0.000000, 0.000000, 0.000000}
\pgfsetstrokecolor{dialinecolor}
\pgfpathellipse{\pgfpoint{27.816377\du}{23.930519\du}}{\pgfpoint{3.725000\du}{0\du}}{\pgfpoint{0\du}{3.825000\du}}
\pgfusepath{stroke}
% setfont left to latex
\definecolor{dialinecolor}{rgb}{0.000000, 0.000000, 0.000000}
\pgfsetstrokecolor{dialinecolor}
\node[anchor=west] at (24\du,24.593019\du){$\mat{U}_2$};
\pgfsetlinewidth{0.200000\du}
\pgfsetdash{}{0pt}
\pgfsetdash{}{0pt}
\pgfsetbuttcap
{
\definecolor{dialinecolor}{rgb}{0.000000, 0.000000, 0.000000}
\pgfsetfillcolor{dialinecolor}
% was here!!!
\definecolor{dialinecolor}{rgb}{0.000000, 0.000000, 0.000000}
\pgfsetstrokecolor{dialinecolor}
\draw (27.816377\du,27.755519\du)--(27.765881\du,32.442013\du);
}
\pgfsetlinewidth{0.100000\du}
\pgfsetdash{}{0pt}
\pgfsetdash{}{0pt}
\pgfsetbuttcap
\pgfsetmiterjoin
\pgfsetlinewidth{0.100000\du}
\pgfsetbuttcap
\pgfsetmiterjoin
\pgfsetdash{}{0pt}
\definecolor{dialinecolor}{rgb}{0.000000, 0.000000, 0.000000}
\pgfsetfillcolor{dialinecolor}
\pgfpathellipse{\pgfpoint{34.902393\du}{23.683396\du}}{\pgfpoint{0.400000\du}{0\du}}{\pgfpoint{0\du}{0.400000\du}}
\pgfusepath{fill}
\definecolor{dialinecolor}{rgb}{0.000000, 0.000000, 0.000000}
\pgfsetstrokecolor{dialinecolor}
\pgfpathellipse{\pgfpoint{34.902393\du}{23.683396\du}}{\pgfpoint{0.400000\du}{0\du}}{\pgfpoint{0\du}{0.400000\du}}
\pgfusepath{stroke}
\pgfsetbuttcap
\pgfsetmiterjoin
\pgfsetdash{}{0pt}
\definecolor{dialinecolor}{rgb}{0.000000, 0.000000, 0.000000}
\pgfsetstrokecolor{dialinecolor}
\pgfpathellipse{\pgfpoint{34.902393\du}{23.683396\du}}{\pgfpoint{0.400000\du}{0\du}}{\pgfpoint{0\du}{0.400000\du}}
\pgfusepath{stroke}
\pgfsetlinewidth{0.100000\du}
\pgfsetdash{}{0pt}
\pgfsetdash{}{0pt}
\pgfsetbuttcap
\pgfsetmiterjoin
\pgfsetlinewidth{0.100000\du}
\pgfsetbuttcap
\pgfsetmiterjoin
\pgfsetdash{}{0pt}
\definecolor{dialinecolor}{rgb}{0.000000, 0.000000, 0.000000}
\pgfsetfillcolor{dialinecolor}
\pgfpathellipse{\pgfpoint{36.547393\du}{23.683396\du}}{\pgfpoint{0.400000\du}{0\du}}{\pgfpoint{0\du}{0.400000\du}}
\pgfusepath{fill}
\definecolor{dialinecolor}{rgb}{0.000000, 0.000000, 0.000000}
\pgfsetstrokecolor{dialinecolor}
\pgfpathellipse{\pgfpoint{36.547393\du}{23.683396\du}}{\pgfpoint{0.400000\du}{0\du}}{\pgfpoint{0\du}{0.400000\du}}
\pgfusepath{stroke}
\pgfsetbuttcap
\pgfsetmiterjoin
\pgfsetdash{}{0pt}
\definecolor{dialinecolor}{rgb}{0.000000, 0.000000, 0.000000}
\pgfsetstrokecolor{dialinecolor}
\pgfpathellipse{\pgfpoint{36.547393\du}{23.683396\du}}{\pgfpoint{0.400000\du}{0\du}}{\pgfpoint{0\du}{0.400000\du}}
\pgfusepath{stroke}
\pgfsetlinewidth{0.100000\du}
\pgfsetdash{}{0pt}
\pgfsetdash{}{0pt}
\pgfsetbuttcap
\pgfsetmiterjoin
\pgfsetlinewidth{0.100000\du}
\pgfsetbuttcap
\pgfsetmiterjoin
\pgfsetdash{}{0pt}
\definecolor{dialinecolor}{rgb}{0.000000, 0.000000, 0.000000}
\pgfsetfillcolor{dialinecolor}
\pgfpathellipse{\pgfpoint{38.142393\du}{23.683396\du}}{\pgfpoint{0.400000\du}{0\du}}{\pgfpoint{0\du}{0.400000\du}}
\pgfusepath{fill}
\definecolor{dialinecolor}{rgb}{0.000000, 0.000000, 0.000000}
\pgfsetstrokecolor{dialinecolor}
\pgfpathellipse{\pgfpoint{38.142393\du}{23.683396\du}}{\pgfpoint{0.400000\du}{0\du}}{\pgfpoint{0\du}{0.400000\du}}
\pgfusepath{stroke}
\pgfsetbuttcap
\pgfsetmiterjoin
\pgfsetdash{}{0pt}
\definecolor{dialinecolor}{rgb}{0.000000, 0.000000, 0.000000}
\pgfsetstrokecolor{dialinecolor}
\pgfpathellipse{\pgfpoint{38.142393\du}{23.683396\du}}{\pgfpoint{0.400000\du}{0\du}}{\pgfpoint{0\du}{0.400000\du}}
\pgfusepath{stroke}
\pgfsetlinewidth{0.200000\du}
\pgfsetdash{}{0pt}
\pgfsetdash{}{0pt}
\pgfsetbuttcap
{
\definecolor{dialinecolor}{rgb}{0.000000, 0.000000, 0.000000}
\pgfsetfillcolor{dialinecolor}
% was here!!!
\definecolor{dialinecolor}{rgb}{0.000000, 0.000000, 0.000000}
\pgfsetstrokecolor{dialinecolor}
\draw (46.122875\du,12.092055\du)--(46.170776\du,20.008450\du);
}
\definecolor{dialinecolor}{rgb}{1.000000, 1.000000, 1.000000}
\pgfsetfillcolor{dialinecolor}
\pgfpathellipse{\pgfpoint{46.194509\du}{23.930519\du}}{\pgfpoint{3.725000\du}{0\du}}{\pgfpoint{0\du}{3.825000\du}}
\pgfusepath{fill}
\pgfsetlinewidth{0.200000\du}
\pgfsetdash{}{0pt}
\pgfsetdash{}{0pt}
\definecolor{dialinecolor}{rgb}{0.000000, 0.000000, 0.000000}
\pgfsetstrokecolor{dialinecolor}
\pgfpathellipse{\pgfpoint{46.194509\du}{23.930519\du}}{\pgfpoint{3.725000\du}{0\du}}{\pgfpoint{0\du}{3.825000\du}}
\pgfusepath{stroke}
% setfont left to latex
\definecolor{dialinecolor}{rgb}{0.000000, 0.000000, 0.000000}
\pgfsetstrokecolor{dialinecolor}
\node[anchor=west] at (42\du,24.593019\du){$\mat{U}_D$};
\pgfsetlinewidth{0.200000\du}
\pgfsetdash{}{0pt}
\pgfsetdash{}{0pt}
\pgfsetbuttcap
{
\definecolor{dialinecolor}{rgb}{0.000000, 0.000000, 0.000000}
\pgfsetfillcolor{dialinecolor}
% was here!!!
\definecolor{dialinecolor}{rgb}{0.000000, 0.000000, 0.000000}
\pgfsetstrokecolor{dialinecolor}
\draw (46.194509\du,27.755519\du)--(46.241769\du,31.953233\du);
}
% setfont left to latex
\definecolor{dialinecolor}{rgb}{0.000000, 0.000000, 0.000000}
\pgfsetstrokecolor{dialinecolor}
\node[anchor=west] at (14\du,18.370414\du){$S_1$};
% setfont left to latex
\definecolor{dialinecolor}{rgb}{0.000000, 0.000000, 0.000000}
\pgfsetstrokecolor{dialinecolor}
\node[anchor=west] at (22\du,18.463477\du){$S_2$};
% setfont left to latex
\definecolor{dialinecolor}{rgb}{0.000000, 0.000000, 0.000000}
\pgfsetstrokecolor{dialinecolor}
\node[anchor=west] at (40.214221\du,18.370414\du){$S_D$};
% setfont left to latex
\definecolor{dialinecolor}{rgb}{0.000000, 0.000000, 0.000000}
\pgfsetstrokecolor{dialinecolor}
\node[anchor=west] at (14\du,30.873409\du){$I_1$};
% setfont left to latex
\definecolor{dialinecolor}{rgb}{0.000000, 0.000000, 0.000000}
\pgfsetstrokecolor{dialinecolor}
\node[anchor=west] at (23\du,30.873409\du){$I_2$};
% setfont left to latex
\definecolor{dialinecolor}{rgb}{0.000000, 0.000000, 0.000000}
\pgfsetstrokecolor{dialinecolor}
\node[anchor=west] at (41\du,30.873409\du){$I_D$};
\end{tikzpicture}
\end{minipage}}
\hspace{-13ex}
\subfigure[Diagram of a Tensor Train.]{
\begin{minipage}[b]{0.7\linewidth}
\centering
\input{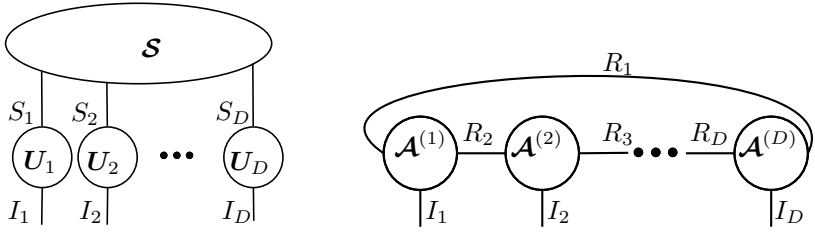}
\end{minipage}}
\caption{Diagram representation of the Tucker and Tensor Train decompositions.}
\label{fig:decomp}
\end{figure}

The Tucker decomposition represents a tensor $\ten{A} \in \mathbb{R}^{I_1\times \cdots \times I_D}$ as
\begin{align}
\label{eqn:Tucker}
\ten{A} &= \ten{S} \times_1 \mat{U}_1 \times_2 \cdots \times_D \mat{U}_D,
\end{align}
where $\ten{S} \in \mathbb{R}^{S_1 \times \cdots \times S_d}$ is called the Tucker core tensor and $\mat{U}_d \in \mathbb{R}^{I_d \times S_d}\,(1\leq d \leq D)$ are the Tucker factor matrices. The total storage complexity of the Tucker decomposition is therefore $\prod_{d=1}^DS_d + \sum_{d=1}^D I_dS_d$. These factor matrices are typically chosen to be orthogonal and can then be obtained as the left singular vectors of the corresponding unfolded matrices of $\ten{A}$. A special case of the Tucker decomposition is the HOSVD~\cite{HOSVD2000}, which has orthogonal matrices and where the Tucker core satisfies two additional properties. The dimensions $S_1,\ldots,S_D$ of the Tucker core are called the multilinear rank of $\ten{A}$ and are defined as
\begin{align*}
    S_d &= \textrm{rank}(\mat{A}_{<d>})\; \leq \; I_d
\end{align*}
for all values of $d$. A graphical representation of the Tucker decomposition is shown in Figure~\ref{fig:decomp}(a).
\subsection{Tensor Train decomposition}
The TT decomposition was introduced into the scientific computing community in~\cite{ivanTT}, but was known as a Matrix Product State in the field of condensed matter physics~\cite{MPS1997,SCHOLLWOCK201196} a decade earlier. 
\begin{definition}
\label{def:TT}
The TT decomposition of a given tensor $\ten{A} \in\mathbb{R}^{I_1\times I_2\times \cdots\times I_D}$ is a set of 3-way tensors $\ten{A}^{(d)} \in \mathbb{R}^{R_d \times I_d \times R_{d+1}}\; (1 \leq d \leq D)$ with $R_1=R_{D+1}=1$ such that each entry $\ten{A}(i_1,i_2,\ldots,i_D)$ can be computed from
\begin{align}
     & \sum_{r_1=1}^{R_1} \sum_{r_2=1}^{R_2} \cdots \sum_{r_D=1}^{R_D} \ten{A}^{(1)}(r_1,i_1,r_2)\,\ten{A}^{(2)}(r_2,i_2,r_3)\, \ldots, \,\ten{A}^{(D)}(r_D,i_D,r_1).
     \label{eq:TTdef}
\end{align}
\end{definition}
The 3-way tensors of the TT are also called the TT-cores and the minimal values of $R_1,\ldots,R_D$ for which \eqref{eq:TTdef} holds exactly for all tensor entries are called the TT-ranks. When $R_1=R_{D+1} > 1$ the decomposition is called a Tensor Ring (TR), for which the diagram is shown in Figure~\ref{fig:decomp}(b) with all dimensions of the TT-cores indicated. We will consider from now on only the TT case and therefore the $R_1$-link in Figure~\ref{fig:decomp}(b) that ``closes the loop" will not be drawn in future diagrams anymore. The total storage complexity of a TT is $\sum_{d=1}^D R_dR_{d+1}I_d$. The TT-ranks are upper bounded as described by the following theorem.
\begin{theorem}(Theorem 2.1 of \cite{ttcross})
\label{theo:TTranks}
For any tensor $\ten{A} \in \mathbb{R}^{I_1 \times \cdots \times I_D}$ there exists a TT-decomposition with TT-ranks
\begin{align*}
R_d \leq \textrm{min}\,\left( \prod_{k=1}^{d-1} I_k, \prod_{k=d}^D I_k \right),
\end{align*}
for $d=2,\ldots,D-1$.
\end{theorem}
Suppose now that we have a Tucker core in TT form. The mode-products of the Tucker factor matrices with this Tucker core in TT form do not alter its TT-ranks. Theorem~\ref{theo:TTranks} therefore reveals the connection between the upper bounds on the TT-ranks of a given tensor and its multilinear rank.
\begin{corollary}
Let $\ten{A}$ be a $D$-way tensor with multilinear rank $S_1,\ldots,S_D$, then its TT-ranks $R_2,\ldots,R_{D}$ satisfy
\begin{align*}
R_d \leq \textrm{min}\,\left( \prod_{k=1}^{d-1} S_k, \prod_{k=d}^D S_k \right),
\end{align*}
for $d=2,\ldots,D-1$.
\end{corollary}

The TT approximation of a given tensor with a prescribed relative error can be computed with either the TT-SVD algorithm~\cite[p.~2301]{ivanTT} or TT-cross algorithm~\cite{ttcross}. Furthermore, through the TT-rounding procedure~\cite[p.~2305]{ivanTT} the TT-ranks of a given TT can be truncated such that the computed approximation satisfies a prescribed relative error. The notion of a TT in site-$d$-mixed-canonical form will be very important in the development of the algorithms in this article and relies on both left-orthogonal and right-orthogonal TT-cores.
 \begin{definition}(~\cite[p.~A689]{Holtz2012})
 A TT-core $\ten{A}^{(d)}$ is left-orthogonal if it can be reshaped into an $R_{d}I_d \times R_{d+1}$ matrix $\mat{A}_d$ such that
 \begin{align*}
 \mat{A}_d^T \, \mat{A}_d &= \mat{I}_{R_d}.
\end{align*}
Similarly, a TT-core $\ten{A}^{(d)}$ is right-orthogonal if it can be reshaped into an $R_{d}\times I_d R_{d+1}$ matrix $\mat{\tilde{A}}_d$ such that
 \begin{align*}
 \mat{\tilde{A}}_d \, \mat{\tilde{A}}_d^T &= \mat{I}_{R_{d-1}}.
\end{align*}
A TT is in site-$d$-mixed-canonical form when all TT-cores $\ten{A}^{(1)}$ up to $\ten{A}^{(d-1)}$ are left-orthogonal and all TT-cores $\ten{A}^{(d+1)}$ up to $\ten{A}^{(D)}$ are right-orthogonal.
\end{definition}
Once a TT is in site-$d$-mixed-canonical form, then it can be readily verified that its Frobenius norm is easily obtained from the $d$th core tensor
\begin{align*}
    ||\ten{A}||_F^2 &= ||\ten{A}^{(d)}||_F^2.
\end{align*}

\subsection{MERA}
\label{subsec:mera}
\begin{figure}[t]
\centering
\subfigure[Binary Hierarchical Tucker.]{
\begin{minipage}[b]{0.45\linewidth}
% Graphic for TeX using PGF
% Title: /home/kbatselier/Pictures/Dia/binaryHT.dia
% Creator: Dia v0.97.3
% CreationDate: Wed Nov 21 13:58:30 2018
% For: kbatselier
% \usepackage{tikz}
% The following commands are not supported in PSTricks at present
% We define them conditionally, so when they are implemented,
% this pgf file will use them.
\ifx\du\undefined
  \newlength{\du}
\fi
\setlength{\du}{4\unitlength}
\begin{tikzpicture}
\pgftransformxscale{1.000000}
\pgftransformyscale{-1.000000}
\definecolor{dialinecolor}{rgb}{0.000000, 0.000000, 0.000000}
\pgfsetstrokecolor{dialinecolor}
\definecolor{dialinecolor}{rgb}{1.000000, 1.000000, 1.000000}
\pgfsetfillcolor{dialinecolor}
\definecolor{dialinecolor}{rgb}{1.000000, 1.000000, 1.000000}
\pgfsetfillcolor{dialinecolor}
\pgfpathellipse{\pgfpoint{47.088250\du}{-157.670381\du}}{\pgfpoint{2.165043\du}{0\du}}{\pgfpoint{0\du}{2.165043\du}}
\pgfusepath{fill}
\pgfsetlinewidth{0.200000\du}
\pgfsetdash{}{0pt}
\pgfsetdash{}{0pt}
\definecolor{dialinecolor}{rgb}{0.000000, 0.000000, 0.000000}
\pgfsetstrokecolor{dialinecolor}
\pgfpathellipse{\pgfpoint{47.088250\du}{-157.670381\du}}{\pgfpoint{2.165043\du}{0\du}}{\pgfpoint{0\du}{2.165043\du}}
\pgfusepath{stroke}
\pgfsetlinewidth{0.200000\du}
\pgfsetdash{}{0pt}
\pgfsetdash{}{0pt}
\pgfsetbuttcap
{
\definecolor{dialinecolor}{rgb}{0.000000, 0.000000, 0.000000}
\pgfsetfillcolor{dialinecolor}
% was here!!!
\definecolor{dialinecolor}{rgb}{0.000000, 0.000000, 0.000000}
\pgfsetstrokecolor{dialinecolor}
\draw (44.923207\du,-157.670381\du)--(39.443263\du,-150.927890\du);
}
\pgfsetlinewidth{0.200000\du}
\pgfsetdash{}{0pt}
\pgfsetdash{}{0pt}
\pgfsetbuttcap
{
\definecolor{dialinecolor}{rgb}{0.000000, 0.000000, 0.000000}
\pgfsetfillcolor{dialinecolor}
% was here!!!
\definecolor{dialinecolor}{rgb}{0.000000, 0.000000, 0.000000}
\pgfsetstrokecolor{dialinecolor}
\draw (54.991970\du,-150.927890\du)--(49.253293\du,-157.670381\du);
}
\definecolor{dialinecolor}{rgb}{1.000000, 1.000000, 1.000000}
\pgfsetfillcolor{dialinecolor}
\pgfpathellipse{\pgfpoint{39.443263\du}{-148.762847\du}}{\pgfpoint{2.165043\du}{0\du}}{\pgfpoint{0\du}{2.165043\du}}
\pgfusepath{fill}
\pgfsetlinewidth{0.200000\du}
\pgfsetdash{}{0pt}
\pgfsetdash{}{0pt}
\definecolor{dialinecolor}{rgb}{0.000000, 0.000000, 0.000000}
\pgfsetstrokecolor{dialinecolor}
\pgfpathellipse{\pgfpoint{39.443263\du}{-148.762847\du}}{\pgfpoint{2.165043\du}{0\du}}{\pgfpoint{0\du}{2.165043\du}}
\pgfusepath{stroke}
\pgfsetlinewidth{0.200000\du}
\pgfsetdash{}{0pt}
\pgfsetdash{}{0pt}
\pgfsetbuttcap
{
\definecolor{dialinecolor}{rgb}{0.000000, 0.000000, 0.000000}
\pgfsetfillcolor{dialinecolor}
% was here!!!
\definecolor{dialinecolor}{rgb}{0.000000, 0.000000, 0.000000}
\pgfsetstrokecolor{dialinecolor}
\draw (37.278220\du,-148.762847\du)--(35.039098\du,-143.282108\du);
}
\pgfsetlinewidth{0.200000\du}
\pgfsetdash{}{0pt}
\pgfsetdash{}{0pt}
\pgfsetbuttcap
{
\definecolor{dialinecolor}{rgb}{0.000000, 0.000000, 0.000000}
\pgfsetfillcolor{dialinecolor}
% was here!!!
\definecolor{dialinecolor}{rgb}{0.000000, 0.000000, 0.000000}
\pgfsetstrokecolor{dialinecolor}
\draw (43.239606\du,-143.282108\du)--(41.608306\du,-148.762847\du);
}
\definecolor{dialinecolor}{rgb}{1.000000, 1.000000, 1.000000}
\pgfsetfillcolor{dialinecolor}
\pgfpathellipse{\pgfpoint{35.039098\du}{-141.117065\du}}{\pgfpoint{2.165043\du}{0\du}}{\pgfpoint{0\du}{2.165043\du}}
\pgfusepath{fill}
\pgfsetlinewidth{0.200000\du}
\pgfsetdash{}{0pt}
\pgfsetdash{}{0pt}
\definecolor{dialinecolor}{rgb}{0.000000, 0.000000, 0.000000}
\pgfsetstrokecolor{dialinecolor}
\pgfpathellipse{\pgfpoint{35.039098\du}{-141.117065\du}}{\pgfpoint{2.165043\du}{0\du}}{\pgfpoint{0\du}{2.165043\du}}
\pgfusepath{stroke}
\pgfsetlinewidth{0.200000\du}
\pgfsetdash{}{0pt}
\pgfsetdash{}{0pt}
\pgfsetbuttcap
{
\definecolor{dialinecolor}{rgb}{0.000000, 0.000000, 0.000000}
\pgfsetfillcolor{dialinecolor}
% was here!!!
\definecolor{dialinecolor}{rgb}{0.000000, 0.000000, 0.000000}
\pgfsetstrokecolor{dialinecolor}
\draw (32.874055\du,-141.117065\du)--(32.121384\du,-137.308264\du);
}
\pgfsetlinewidth{0.200000\du}
\pgfsetdash{}{0pt}
\pgfsetdash{}{0pt}
\pgfsetbuttcap
{
\definecolor{dialinecolor}{rgb}{0.000000, 0.000000, 0.000000}
\pgfsetfillcolor{dialinecolor}
% was here!!!
\definecolor{dialinecolor}{rgb}{0.000000, 0.000000, 0.000000}
\pgfsetstrokecolor{dialinecolor}
\draw (37.961012\du,-137.356928\du)--(37.204141\du,-141.117065\du);
}
\definecolor{dialinecolor}{rgb}{1.000000, 1.000000, 1.000000}
\pgfsetfillcolor{dialinecolor}
\pgfpathellipse{\pgfpoint{43.239606\du}{-141.117065\du}}{\pgfpoint{2.165043\du}{0\du}}{\pgfpoint{0\du}{2.165043\du}}
\pgfusepath{fill}
\pgfsetlinewidth{0.200000\du}
\pgfsetdash{}{0pt}
\pgfsetdash{}{0pt}
\definecolor{dialinecolor}{rgb}{0.000000, 0.000000, 0.000000}
\pgfsetstrokecolor{dialinecolor}
\pgfpathellipse{\pgfpoint{43.239606\du}{-141.117065\du}}{\pgfpoint{2.165043\du}{0\du}}{\pgfpoint{0\du}{2.165043\du}}
\pgfusepath{stroke}
\pgfsetlinewidth{0.200000\du}
\pgfsetdash{}{0pt}
\pgfsetdash{}{0pt}
\pgfsetbuttcap
{
\definecolor{dialinecolor}{rgb}{0.000000, 0.000000, 0.000000}
\pgfsetfillcolor{dialinecolor}
% was here!!!
\definecolor{dialinecolor}{rgb}{0.000000, 0.000000, 0.000000}
\pgfsetstrokecolor{dialinecolor}
\draw (41.074563\du,-141.117065\du)--(40.321892\du,-137.308264\du);
}
\pgfsetlinewidth{0.200000\du}
\pgfsetdash{}{0pt}
\pgfsetdash{}{0pt}
\pgfsetbuttcap
{
\definecolor{dialinecolor}{rgb}{0.000000, 0.000000, 0.000000}
\pgfsetfillcolor{dialinecolor}
% was here!!!
\definecolor{dialinecolor}{rgb}{0.000000, 0.000000, 0.000000}
\pgfsetstrokecolor{dialinecolor}
\draw (46.161521\du,-137.356928\du)--(45.404649\du,-141.117065\du);
}
\definecolor{dialinecolor}{rgb}{1.000000, 1.000000, 1.000000}
\pgfsetfillcolor{dialinecolor}
\pgfpathellipse{\pgfpoint{54.991970\du}{-148.762847\du}}{\pgfpoint{2.165043\du}{0\du}}{\pgfpoint{0\du}{2.165043\du}}
\pgfusepath{fill}
\pgfsetlinewidth{0.200000\du}
\pgfsetdash{}{0pt}
\pgfsetdash{}{0pt}
\definecolor{dialinecolor}{rgb}{0.000000, 0.000000, 0.000000}
\pgfsetstrokecolor{dialinecolor}
\pgfpathellipse{\pgfpoint{54.991970\du}{-148.762847\du}}{\pgfpoint{2.165043\du}{0\du}}{\pgfpoint{0\du}{2.165043\du}}
\pgfusepath{stroke}
\pgfsetlinewidth{0.200000\du}
\pgfsetdash{}{0pt}
\pgfsetdash{}{0pt}
\pgfsetbuttcap
{
\definecolor{dialinecolor}{rgb}{0.000000, 0.000000, 0.000000}
\pgfsetfillcolor{dialinecolor}
% was here!!!
\definecolor{dialinecolor}{rgb}{0.000000, 0.000000, 0.000000}
\pgfsetstrokecolor{dialinecolor}
\draw (52.826927\du,-148.762847\du)--(50.587805\du,-143.282108\du);
}
\pgfsetlinewidth{0.200000\du}
\pgfsetdash{}{0pt}
\pgfsetdash{}{0pt}
\pgfsetbuttcap
{
\definecolor{dialinecolor}{rgb}{0.000000, 0.000000, 0.000000}
\pgfsetfillcolor{dialinecolor}
% was here!!!
\definecolor{dialinecolor}{rgb}{0.000000, 0.000000, 0.000000}
\pgfsetstrokecolor{dialinecolor}
\draw (58.788313\du,-143.282108\du)--(57.157013\du,-148.762847\du);
}
\definecolor{dialinecolor}{rgb}{1.000000, 1.000000, 1.000000}
\pgfsetfillcolor{dialinecolor}
\pgfpathellipse{\pgfpoint{50.587805\du}{-141.117065\du}}{\pgfpoint{2.165043\du}{0\du}}{\pgfpoint{0\du}{2.165043\du}}
\pgfusepath{fill}
\pgfsetlinewidth{0.200000\du}
\pgfsetdash{}{0pt}
\pgfsetdash{}{0pt}
\definecolor{dialinecolor}{rgb}{0.000000, 0.000000, 0.000000}
\pgfsetstrokecolor{dialinecolor}
\pgfpathellipse{\pgfpoint{50.587805\du}{-141.117065\du}}{\pgfpoint{2.165043\du}{0\du}}{\pgfpoint{0\du}{2.165043\du}}
\pgfusepath{stroke}
\pgfsetlinewidth{0.200000\du}
\pgfsetdash{}{0pt}
\pgfsetdash{}{0pt}
\pgfsetbuttcap
{
\definecolor{dialinecolor}{rgb}{0.000000, 0.000000, 0.000000}
\pgfsetfillcolor{dialinecolor}
% was here!!!
\definecolor{dialinecolor}{rgb}{0.000000, 0.000000, 0.000000}
\pgfsetstrokecolor{dialinecolor}
\draw (48.422761\du,-141.117065\du)--(47.670091\du,-137.308264\du);
}
\pgfsetlinewidth{0.200000\du}
\pgfsetdash{}{0pt}
\pgfsetdash{}{0pt}
\pgfsetbuttcap
{
\definecolor{dialinecolor}{rgb}{0.000000, 0.000000, 0.000000}
\pgfsetfillcolor{dialinecolor}
% was here!!!
\definecolor{dialinecolor}{rgb}{0.000000, 0.000000, 0.000000}
\pgfsetstrokecolor{dialinecolor}
\draw (53.509719\du,-137.356928\du)--(52.752848\du,-141.117065\du);
}
\definecolor{dialinecolor}{rgb}{1.000000, 1.000000, 1.000000}
\pgfsetfillcolor{dialinecolor}
\pgfpathellipse{\pgfpoint{58.788313\du}{-141.117065\du}}{\pgfpoint{2.165043\du}{0\du}}{\pgfpoint{0\du}{2.165043\du}}
\pgfusepath{fill}
\pgfsetlinewidth{0.200000\du}
\pgfsetdash{}{0pt}
\pgfsetdash{}{0pt}
\definecolor{dialinecolor}{rgb}{0.000000, 0.000000, 0.000000}
\pgfsetstrokecolor{dialinecolor}
\pgfpathellipse{\pgfpoint{58.788313\du}{-141.117065\du}}{\pgfpoint{2.165043\du}{0\du}}{\pgfpoint{0\du}{2.165043\du}}
\pgfusepath{stroke}
\pgfsetlinewidth{0.200000\du}
\pgfsetdash{}{0pt}
\pgfsetdash{}{0pt}
\pgfsetbuttcap
{
\definecolor{dialinecolor}{rgb}{0.000000, 0.000000, 0.000000}
\pgfsetfillcolor{dialinecolor}
% was here!!!
\definecolor{dialinecolor}{rgb}{0.000000, 0.000000, 0.000000}
\pgfsetstrokecolor{dialinecolor}
\draw (56.623269\du,-141.117065\du)--(55.870599\du,-137.308264\du);
}
\pgfsetlinewidth{0.200000\du}
\pgfsetdash{}{0pt}
\pgfsetdash{}{0pt}
\pgfsetbuttcap
{
\definecolor{dialinecolor}{rgb}{0.000000, 0.000000, 0.000000}
\pgfsetfillcolor{dialinecolor}
% was here!!!
\definecolor{dialinecolor}{rgb}{0.000000, 0.000000, 0.000000}
\pgfsetstrokecolor{dialinecolor}
\draw (61.710227\du,-137.356928\du)--(60.953356\du,-141.117065\du);
}
\end{tikzpicture}
\end{minipage}}
\subfigure[Binary MERA.]{
\begin{minipage}[b]{0.45\linewidth}
\input{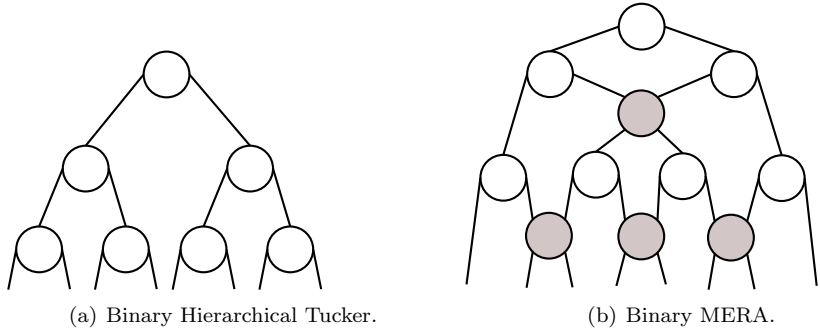}
\end{minipage}}
\caption{Diagram representation of two hierarchical tensor decompositions.}
\label{fig:decomp2}
\end{figure}

\begin{figure}[t]
    \centering
    \input{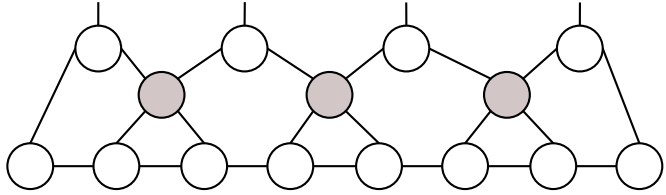}
    \caption{A TT of an 8-way tensor (bottom row) is coarse-grained into a 4-way tensor through one layer of a MERA.}
    \label{fig:coarsegrainingMERA}
\end{figure}

The MERA decomposition is a generalization of the Hierarchical Tucker decomposition and consists of three different building blocks. A common implementation of the Hierarchical Tucker decomposition is the binary tree form, as shown in Figure~\ref{fig:decomp2}(a). Reading such a diagram from the bottom to the top, one can interpret each row/layer in such a tree structure as a coarse-graining transformation where each tensor in a row/layer transforms two indices into one index. Such tensors $\ten{W}$ of size $I_1 \times \cdots \times I_K \times S$ that reduce $K>1$ indices to a single index are called isometries. An isometry can always be reshaped into a size $I_1I_2\cdots I_K \times S$ matrix $\mat{W}$ with orthonormal columns
\begin{align*}
    \mat{W}^T\;\mat{W} &= \mat{I}_S,
\end{align*}
where $S$ is the dimension of the ``output" index. The minimal outgoing dimensions of all isometries such that the MERA represents a given tensor exactly are called the MERA-ranks. The diagram representation of an isometry is shown in Figure~\ref{fig:MERAunits}(a). The bottom layer of isometries with $K=2$ in Figure~\ref{fig:decomp2}(a) reduces the eight indices of a given tensor into four indices, as illustrated in Figure~\ref{fig:coarsegrainingHT}. Each application of a layer in the tree halves the resulting total number of indices. The coarse-graining with a Hierarchical Tucker decomposition pairs two consecutive indices and sums over them, thereby ignoring possible correlations over neighbouring indices resulting in higher ranks during coarse-graining. This issue is resolved in the MERA through the introduction of additional disentangler tensors in the coarse-graining layers. Disentanglers, shown as shaded nodes in Figure~\ref{fig:decomp2}(b), ``bridge" neighbouring pairs before being coarse-grained. A disentangler tensor is per definition a 4-way tensor $\ten{V}$ of size $I_1 \times I_2 \times I_1 \times I_2$ that can be reshaped into an orthogonal $I_1I_2 \times I_1I_2$ matrix $\mat{V}$. The reduction of an 8-way TT into a 4-way TT through a MERA layer is shown as a diagram in Figure~\ref{fig:coarsegrainingMERA}. The third and final MERA building block is the top tensor. This tensor $\ten{T}$ is located at the top of the MERA structure and connects to all outgoing isometry indices of the highest layer. Since all disentanglers and isometries have their respective notion of orthogonality, it follows that the Frobenius norm of a tensor $\ten{A}$ that is represented by a MERA is given by $||\ten{A}||_F^2=||\ten{T}||_F^2$. This easy computation of the norm due to orthogonality is very similar to the case of a TT in site-$d$-mixed canonical form. The storage complexity of a MERA is simply the sum of storage complexities of all disentanglers, isometries and the top tensor. In this respect, it is only meaningful from a data tensor compression perspective to have MERA-ranks that do not increase over consecutive layers. In the next section, we develop the main algorithm to convert a given TT into a Tucker decomposition and this algorithm will serve as the main computational building block to eventually convert a TT into a MERA.

\begin{figure}[t]
\subfigure[Isometry tensor.]{
\begin{minipage}[b]{0.45\linewidth}
\centering
% Graphic for TeX using PGF
% Title: /home/kbatselier/Pictures/Dia/Isometry.dia
% Creator: Dia v0.97.3
% CreationDate: Tue Dec  4 10:23:43 2018
% For: kbatselier
% \usepackage{tikz}
% The following commands are not supported in PSTricks at present
% We define them conditionally, so when they are implemented,
% this pgf file will use them.
\ifx\du\undefined
  \newlength{\du}
\fi
\setlength{\du}{5\unitlength}
\begin{tikzpicture}
\pgftransformxscale{1.000000}
\pgftransformyscale{-1.000000}
\definecolor{dialinecolor}{rgb}{0.000000, 0.000000, 0.000000}
\pgfsetstrokecolor{dialinecolor}
\definecolor{dialinecolor}{rgb}{1.000000, 1.000000, 1.000000}
\pgfsetfillcolor{dialinecolor}
\definecolor{dialinecolor}{rgb}{1.000000, 1.000000, 1.000000}
\pgfsetfillcolor{dialinecolor}
\pgfpathellipse{\pgfpoint{35.492274\du}{9.415089\du}}{\pgfpoint{2.165043\du}{0\du}}{\pgfpoint{0\du}{2.165043\du}}
\pgfusepath{fill}
\pgfsetlinewidth{0.200000\du}
\pgfsetdash{}{0pt}
\pgfsetdash{}{0pt}
\definecolor{dialinecolor}{rgb}{0.000000, 0.000000, 0.000000}
\pgfsetstrokecolor{dialinecolor}
\pgfpathellipse{\pgfpoint{35.492274\du}{9.415089\du}}{\pgfpoint{2.165043\du}{0\du}}{\pgfpoint{0\du}{2.165043\du}}
\pgfusepath{stroke}
\pgfsetlinewidth{0.200000\du}
\pgfsetdash{}{0pt}
\pgfsetdash{}{0pt}
\pgfsetbuttcap
{
\definecolor{dialinecolor}{rgb}{0.000000, 0.000000, 0.000000}
\pgfsetfillcolor{dialinecolor}
% was here!!!
\definecolor{dialinecolor}{rgb}{0.000000, 0.000000, 0.000000}
\pgfsetstrokecolor{dialinecolor}
\draw (33.327231\du,9.415089\du)--(32.350000\du,15.150000\du);
}
\pgfsetlinewidth{0.200000\du}
\pgfsetdash{}{0pt}
\pgfsetdash{}{0pt}
\pgfsetbuttcap
{
\definecolor{dialinecolor}{rgb}{0.000000, 0.000000, 0.000000}
\pgfsetfillcolor{dialinecolor}
% was here!!!
\definecolor{dialinecolor}{rgb}{0.000000, 0.000000, 0.000000}
\pgfsetstrokecolor{dialinecolor}
\draw (38.800000\du,15.150000\du)--(37.657317\du,9.415089\du);
}
\pgfsetlinewidth{0.200000\du}
\pgfsetdash{}{0pt}
\pgfsetdash{}{0pt}
\pgfsetbuttcap
{
\definecolor{dialinecolor}{rgb}{0.000000, 0.000000, 0.000000}
\pgfsetfillcolor{dialinecolor}
% was here!!!
\definecolor{dialinecolor}{rgb}{0.000000, 0.000000, 0.000000}
\pgfsetstrokecolor{dialinecolor}
\draw (35.494723\du,3.928803\du)--(35.492274\du,7.250046\du);
}
\pgfsetlinewidth{0.200000\du}
\pgfsetdash{}{0pt}
\pgfsetdash{}{0pt}
\pgfsetbuttcap
{
\definecolor{dialinecolor}{rgb}{0.000000, 0.000000, 0.000000}
\pgfsetfillcolor{dialinecolor}
% was here!!!
\definecolor{dialinecolor}{rgb}{0.000000, 0.000000, 0.000000}
\pgfsetstrokecolor{dialinecolor}
\draw (33.961358\du,10.946006\du)--(33.957260\du,15.070885\du);
}
\pgfsetlinewidth{0.200000\du}
\pgfsetdash{}{0pt}
\pgfsetdash{}{0pt}
\pgfsetbuttcap
{
\definecolor{dialinecolor}{rgb}{0.000000, 0.000000, 0.000000}
\pgfsetfillcolor{dialinecolor}
% was here!!!
\definecolor{dialinecolor}{rgb}{0.000000, 0.000000, 0.000000}
\pgfsetstrokecolor{dialinecolor}
\draw (37.023191\du,10.946006\du)--(37.032186\du,15.143236\du);
}
% setfont left to latex
\definecolor{dialinecolor}{rgb}{0.000000, 0.000000, 0.000000}
\pgfsetstrokecolor{dialinecolor}
\node[anchor=west] at (33.875000\du,2.715000\du){$S$};
\definecolor{dialinecolor}{rgb}{0.000000, 0.000000, 0.000000}
\pgfsetfillcolor{dialinecolor}
\pgfpathellipse{\pgfpoint{33.994269\du}{16.511511\du}}{\pgfpoint{0.292314\du}{0\du}}{\pgfpoint{0\du}{0.292314\du}}
\pgfusepath{fill}
\pgfsetlinewidth{0.100000\du}
\pgfsetdash{}{0pt}
\pgfsetdash{}{0pt}
\definecolor{dialinecolor}{rgb}{0.000000, 0.000000, 0.000000}
\pgfsetstrokecolor{dialinecolor}
\pgfpathellipse{\pgfpoint{33.994269\du}{16.511511\du}}{\pgfpoint{0.292314\du}{0\du}}{\pgfpoint{0\du}{0.292314\du}}
\pgfusepath{stroke}
\definecolor{dialinecolor}{rgb}{0.000000, 0.000000, 0.000000}
\pgfsetfillcolor{dialinecolor}
\pgfpathellipse{\pgfpoint{35.008678\du}{16.511511\du}}{\pgfpoint{0.292314\du}{0\du}}{\pgfpoint{0\du}{0.292314\du}}
\pgfusepath{fill}
\pgfsetlinewidth{0.100000\du}
\pgfsetdash{}{0pt}
\pgfsetdash{}{0pt}
\definecolor{dialinecolor}{rgb}{0.000000, 0.000000, 0.000000}
\pgfsetstrokecolor{dialinecolor}
\pgfpathellipse{\pgfpoint{35.008678\du}{16.511511\du}}{\pgfpoint{0.292314\du}{0\du}}{\pgfpoint{0\du}{0.292314\du}}
\pgfusepath{stroke}
\definecolor{dialinecolor}{rgb}{0.000000, 0.000000, 0.000000}
\pgfsetfillcolor{dialinecolor}
\pgfpathellipse{\pgfpoint{36.004353\du}{16.511511\du}}{\pgfpoint{0.292314\du}{0\du}}{\pgfpoint{0\du}{0.292314\du}}
\pgfusepath{fill}
\pgfsetlinewidth{0.100000\du}
\pgfsetdash{}{0pt}
\pgfsetdash{}{0pt}
\definecolor{dialinecolor}{rgb}{0.000000, 0.000000, 0.000000}
\pgfsetstrokecolor{dialinecolor}
\pgfpathellipse{\pgfpoint{36.004353\du}{16.511511\du}}{\pgfpoint{0.292314\du}{0\du}}{\pgfpoint{0\du}{0.292314\du}}
\pgfusepath{stroke}
% setfont left to latex
\definecolor{dialinecolor}{rgb}{0.000000, 0.000000, 0.000000}
\pgfsetstrokecolor{dialinecolor}
\node[anchor=west] at (30.325000\du,16.315000\du){$I_1$};
% setfont left to latex
\definecolor{dialinecolor}{rgb}{0.000000, 0.000000, 0.000000}
\pgfsetstrokecolor{dialinecolor}
\node[anchor=west] at (37.850000\du,16.315000\du){$I_K$};
\end{tikzpicture}
\end{minipage}}
\subfigure[Disentangler tensor.]{
\begin{minipage}[b]{0.45\linewidth}
\centering
% Graphic for TeX using PGF
% Title: /home/kbatselier/Pictures/Dia/Disentangler.dia
% Creator: Dia v0.97.3
% CreationDate: Tue Dec  4 10:23:30 2018
% For: kbatselier
% \usepackage{tikz}
% The following commands are not supported in PSTricks at present
% We define them conditionally, so when they are implemented,
% this pgf file will use them.
\ifx\du\undefined
  \newlength{\du}
\fi
\setlength{\du}{5\unitlength}
\begin{tikzpicture}
\pgftransformxscale{1.000000}
\pgftransformyscale{-1.000000}
\definecolor{dialinecolor}{rgb}{0.000000, 0.000000, 0.000000}
\pgfsetstrokecolor{dialinecolor}
\definecolor{dialinecolor}{rgb}{1.000000, 1.000000, 1.000000}
\pgfsetfillcolor{dialinecolor}
\pgfsetlinewidth{0.200000\du}
\pgfsetdash{}{0pt}
\pgfsetdash{}{0pt}
\pgfsetbuttcap
{
\definecolor{dialinecolor}{rgb}{0.000000, 0.000000, 0.000000}
\pgfsetfillcolor{dialinecolor}
% was here!!!
\definecolor{dialinecolor}{rgb}{0.000000, 0.000000, 0.000000}
\pgfsetstrokecolor{dialinecolor}
\draw (15.359126\du,9.008458\du)--(14.742700\du,5.056131\du);
}
\pgfsetlinewidth{0.200000\du}
\pgfsetdash{}{0pt}
\pgfsetdash{}{0pt}
\pgfsetbuttcap
{
\definecolor{dialinecolor}{rgb}{0.000000, 0.000000, 0.000000}
\pgfsetfillcolor{dialinecolor}
% was here!!!
\definecolor{dialinecolor}{rgb}{0.000000, 0.000000, 0.000000}
\pgfsetstrokecolor{dialinecolor}
\draw (19.157600\du,4.765631\du)--(18.420960\du,9.008458\du);
}
\definecolor{dialinecolor}{rgb}{0.823529, 0.776471, 0.776471}
\pgfsetfillcolor{dialinecolor}
\pgfpathellipse{\pgfpoint{16.890043\du}{10.539374\du}}{\pgfpoint{2.165043\du}{0\du}}{\pgfpoint{0\du}{2.165043\du}}
\pgfusepath{fill}
\pgfsetlinewidth{0.200000\du}
\pgfsetdash{}{0pt}
\pgfsetdash{}{0pt}
\definecolor{dialinecolor}{rgb}{0.000000, 0.000000, 0.000000}
\pgfsetstrokecolor{dialinecolor}
\pgfpathellipse{\pgfpoint{16.890043\du}{10.539374\du}}{\pgfpoint{2.165043\du}{0\du}}{\pgfpoint{0\du}{2.165043\du}}
\pgfusepath{stroke}
\pgfsetlinewidth{0.200000\du}
\pgfsetdash{}{0pt}
\pgfsetdash{}{0pt}
\pgfsetbuttcap
{
\definecolor{dialinecolor}{rgb}{0.000000, 0.000000, 0.000000}
\pgfsetfillcolor{dialinecolor}
% was here!!!
\definecolor{dialinecolor}{rgb}{0.000000, 0.000000, 0.000000}
\pgfsetstrokecolor{dialinecolor}
\draw (15.359126\du,12.070291\du)--(14.775000\du,15.464231\du);
}
\pgfsetlinewidth{0.200000\du}
\pgfsetdash{}{0pt}
\pgfsetdash{}{0pt}
\pgfsetbuttcap
{
\definecolor{dialinecolor}{rgb}{0.000000, 0.000000, 0.000000}
\pgfsetfillcolor{dialinecolor}
% was here!!!
\definecolor{dialinecolor}{rgb}{0.000000, 0.000000, 0.000000}
\pgfsetstrokecolor{dialinecolor}
\draw (19.123100\du,15.426331\du)--(18.420960\du,12.070291\du);
}
% setfont left to latex
\definecolor{dialinecolor}{rgb}{0.000000, 0.000000, 0.000000}
\pgfsetstrokecolor{dialinecolor}
\node[anchor=west] at (13.450000\du,17.350000\du){$I_1$};
% setfont left to latex
\definecolor{dialinecolor}{rgb}{0.000000, 0.000000, 0.000000}
\pgfsetstrokecolor{dialinecolor}
\node[anchor=west] at (17.875000\du,17.350000\du){$I_2$};
% setfont left to latex
\definecolor{dialinecolor}{rgb}{0.000000, 0.000000, 0.000000}
\pgfsetstrokecolor{dialinecolor}
\node[anchor=west] at (13.450000\du,3.165000\du){$I_1$};
% setfont left to latex
\definecolor{dialinecolor}{rgb}{0.000000, 0.000000, 0.000000}
\pgfsetstrokecolor{dialinecolor}
\node[anchor=west] at (17.875000\du,3.165000\du){$I_2$};
\end{tikzpicture}
\end{minipage}}
\caption{Diagram representation of two MERA building block tensors.}
\label{fig:MERAunits}
\end{figure}
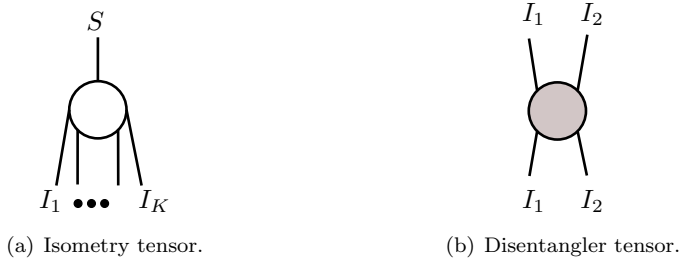

\begin{figure}[t]
    \centering
    \input{figs/coarsegrainingbinaryHT.tex}
    \caption{A TT of an 8-way tensor (bottom row) is coarse-grained into a 4-way tensor through one layer of a HT/TTN.}
    \label{fig:coarsegrainingHT}
\end{figure}

%\begin{figure}[t]
%\centering
%\subfigure[Orthogonality disentanglers.]{
%\begin{minipage}[b]{0.3\linewidth}
%\input{figs/OrthoDisentanglers.tex}
%\end{minipage}}
%\subfigure[Orthogonality isometry.]{
%\begin{minipage}[b]{0.3\linewidth}
%\input{figs/OrthoIsometry1.tex}
%\end{minipage}}
%\caption{Orthogonality isometry.}
%\subfigure[Orthogonality isometry2.]{
%\begin{minipage}[b]{0.3\linewidth}
%\input{figs/OrthoIsometry2.tex}
%\end{minipage}}
%\caption{Orthogonality isometry.}
%\label{fig:ortho}
%\end{figure}

\section{Tensor train to Tucker decomposition}
\label{sec:TT2Tucker}
In this section, an algorithm is developed that converts a given TT into either a HOSVD or truncated HOSVD with a guaranteed upper bound on the relative approximation error. The Tucker core $\ten{S}$ will be directly obtained in the TT format, avoiding its exponential storage complexity. The starting point of the algorithm is a TT in site-$1$-mixed-canonical form. Before stating the algorithm, we first introduce some additional notation together with an important lemma.
\subsection{Tucker factor matrix from TT-core}
In order to know how a given TT can be converted into a Tucker decomposition we need to know how the Tucker factor matrices can be computed from each TT-core. In order to describe this computation we first introduce the following convenient notation.
\begin{definition}
\label{def:TTmatrices}
Let $\ten{A}^{(1)},\ldots,\ten{A}^{(D)}$ be TT-cores of a $D$-way tensor $\ten{A}$. We define $\mat{A}_{<d}$ as the $R_d \times (I_1\cdots I_{d-1})$ matrix obtained from summing over the auxiliary indices of $\ten{A}^{(1)}$ up to $\ten{A}^{(d-1)}$ and permuting and reshaping the result into the desired matrix. The $R_{d+1} \times (I_{d+1}\cdots I_D)$ matrix $\mat{A}_{>d}$ is defined similarly from the TT-cores $\ten{A}^{(d+1)}$ up to $\ten{A}^{(D)}$. The $I_d \times (R_dR_{d+1})$ matrix $\mat{A}_d$ is defined from permuting and reshaping $\ten{A}^{(d)}$. Finally, both $\mat{A}_{<1}$ and $\mat{A}_{>D}$ are defined to be unit scalars.
\end{definition}
Note that if the TT of $\ten{A}$ is in site-$d$-mixed-canonical form, then the left and right-orthogonality of the TT-cores implies that both $\mat{A}_{<d}$ and $\mat{A}_{>d}$ have orthonormal rows
\begin{align*}
    \mat{A}_{<d} \; \mat{A}_{<d}^T = \mat{I}_{R_d}\; \textrm{ and }\; \mat{A}_{>d} \; \mat{A}_{>d}^T &= \mat{I}_{R_{d+1}}.
\end{align*}
The following lemma tells us how the unfolding matrix $\mat{A}_{<d>}$ can be written in terms of the matrices from Definition~\ref{def:TTmatrices}.
\begin{lemma}
\label{lemma:unfolding}
For a $D$-way tensor $\ten{A}$ in TT-form we have the following relationship
\begin{align*}
\mat{A}_{<d>} &= \mat{A}_d \; \left(\mat{A}_{>d} \otimes \mat{A}_{<d} \right) \quad (d=1,\ldots,D).
\end{align*}
\end{lemma}
The Kronecker product in Lemma~\ref{lemma:unfolding} is due to the rank-1 link of the TT. Also note that the rows of $\left(\mat{A}_{>d} \otimes \mat{A}_{<d} \right)$ are orthonormal when the TT is in site-$d$-mixed-canonical form, due to the preservation of orthonormality with the Kronecker product. Lemma~\ref{lemma:unfolding} tells us that any unfolding matrix $\mat{A}_{<d>}$ can be written as a product of $\mat{A}_d$ with $\left(\mat{A}_{>d} \otimes \mat{A}_{<d} \right)$, which leads to the following two corollaries.
\begin{corollary}
\label{cor:multilinearrankbound}
For a $D$-way tensor $\ten{A}$ with multilinear ranks $S_1,\ldots,S_D$ we have that
\begin{align*}
S_d = \textrm{rank}(\mat{A}_{<d>}) = \textrm{rank}(\mat{A}_d) \leq  \textrm{min}(I_d,R_dR_{d+1}) \leq I_d \quad (d=1,\ldots,D).
\end{align*}
\end{corollary}
\begin{corollary}
\label{cor:factormatrix}
For a tensor $\ten{A}$ in site-$d$-mixed-canonical TT-form, let the compact SVD of $\mat{A}_{d}$ be given by $\mat{U}_d\,\mat{S}\,\mat{V}^T$, then the compact SVD of the unfolding matrix $\mat{A}_{<d>}$ is
\begin{align*}
\mat{A}_{<d>} &= \mat{U}_d \; \mat{S} \; \mat{V}^T\left(\mat{A}_{>d} \otimes \mat{A}_{<d} \right).
\end{align*}
\end{corollary}
In Corollary~\ref{cor:multilinearrankbound} we have tacitly assumed that $R_dR_{d+1} < \prod_{k \neq d} I_k$ is always satisfied. Corollary~\ref{cor:factormatrix} follows directly from the fact that the product of matrices with orthonormal rows also has orthonormal rows. The matrix $\left(\mat{A}_{>d}^T \otimes \mat{A}_{<d}^T \right) \mat{V}$ therefore contains the right singular vectors of $\mat{A}_{<d>}$ corresponding with the $S_d$ largest singular values. Corollary~\ref{cor:factormatrix} also implies that the HOSVD factor matrix $\mat{U}_d$ can be directly computed from the SVD of $\mat{A}_d$. The $d$th component $S_d$ of the multilinear rank can be determined by inspecting the singular values on the diagonal $\mat{S}$ matrix. If there is no need to know the exact multilinear rank, then a square orthogonal $\mat{U}_d$ can also be obtained through a QR decomposition of $\mat{A}_d$.
\subsection{The TT to Tucker conversion algorithm}
\begin{figure}[tbh]
    \centering
    \input{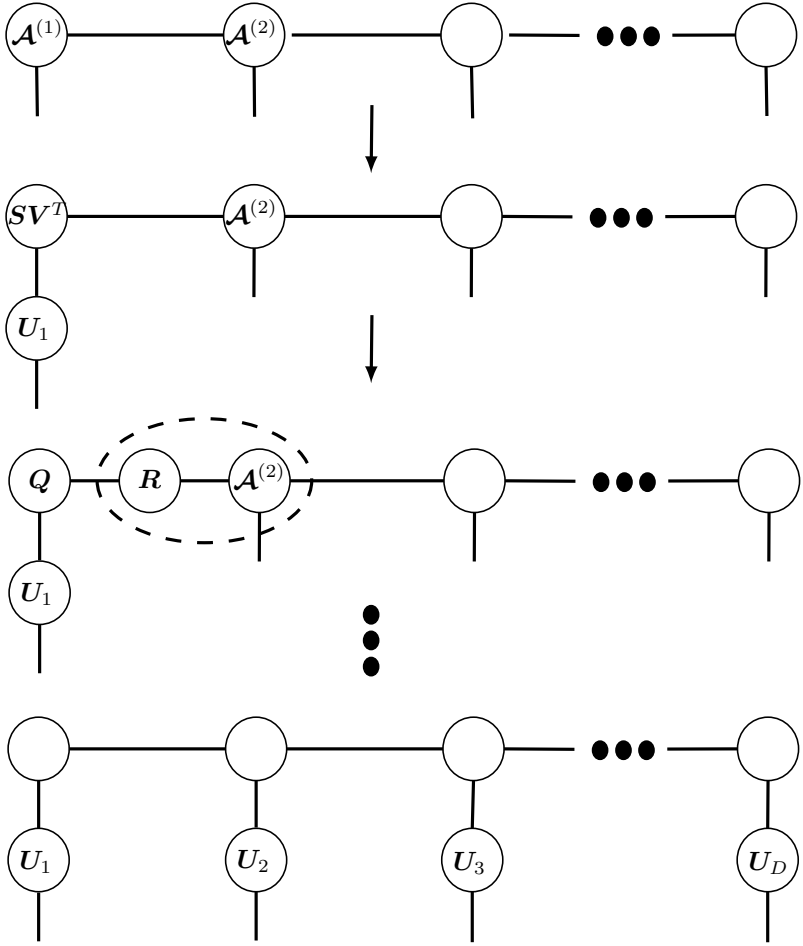}
    \caption{The complete first execution of the for-loop in Algorithm~\ref{alg:TT2Tucker} in diagram form.}
    \label{fig:TT2HOSVD}
\end{figure}

Lemma~\ref{lemma:unfolding} forms the basis of the proposed algorithm to convert a given TT into either a HOSVD or truncated HOSVD. The algorithm to compute a truncated HOSVD is presented in pseudo-code as Algorithm~\ref{alg:TT2Tucker}. The algorithm assumes the TT is in site-1-mixed-canonical form but can be easily adjusted to work for any other starting site. The main idea of Algorithm~\ref{alg:TT2Tucker} is to compute the orthogonal factor matrix $\mat{U}_d$ using Corollary~\ref{cor:factormatrix} and then to bring the TT into site-$(d+1)$-mixed-canonical form.
%The computational complexity for computing all $d$ factor matrices is therefore $\approx O(dR^2N^2)$, where $R$ and $N$ are the maximal TT-rank and dimension of $\ten{A}$, respectively.
The conversion of the TT into site-$d+1$-mixed-canonical form is computed through a QR decomposition of the $\mat{S}\mat{V}^T$ factor.
%, which requires approximately $O(NS^3)$ flops.
The orthogonal $\mat{Q}$ matrix is then retained as the $d$th TT-core of the Tucker core $\ten{S}$, while the norm of $\ten{A}$ is moved to the next TT-core $\ten{A}^{(d+1)}$ through the absorption of the $\mat{R}$ factor. The final TT of $\ten{S}$ will therefore be in site-$D$-mixed-canonical form. Both the SVD step and the QR decomposition step are graphically represented in Figure~\ref{fig:TT2HOSVD}. During each run of the for-loop in Algorithm~\ref{alg:TT2Tucker} we are working with a partially truncated core tensor, which is very reminiscent of the ST-HOSVD algorithm~\cite{STHOSVD2012}. In fact, the approximation error induced by truncating the SVD in Algorithm~\ref{alg:TT2Tucker} can also be expressed exactly in terms of the singular values.
\begin{theorem}
\label{theo:approxerror}
Let $\sigma_d(i_d)$ be the $i_d$th singular value of $\mat{A}_d$ and $\ten{\hat{A}}$ be the tensor computed by Algorithm~\ref{alg:TT2Tucker} with truncated SVDs, then
\begin{align*}
||\ten{A}-\ten{\hat{A}}||_F^2 &= \sum_{d=1}^D \sum_{i_d=S_d+1}^{I_d} \sigma_d(i_d)^2.    
\end{align*}
\end{theorem}
Given that the TT for the all-orthogonal Tucker core is in site-$d$-mixed-canonical form and the similarity of Algorithm~\ref{alg:TT2Tucker} concerning the use of a sequentially truncated Tucker core, it follows that the proof of Theorem~\ref{theo:approxerror} is completely identical with the one found in~\cite[p.~A1039]{STHOSVD2012}. Theorem~\ref{theo:approxerror} allows us to compute the absolute approximation error for a truncated HOSVD during the execution of Algorithm~\ref{alg:TT2Tucker} by simply adding the squares of the discarded singular values. In addition,  Theorem~\ref{theo:approxerror} also allows us to compute a truncated HOSVD for a given upper bound $\epsilon$ on the relative approximation error. Since Algorithm~\ref{alg:TT2Tucker} consists of $D$ truncated SVDs, setting the tolerance $\delta$ for each of these SVDs to $\epsilon ||\ten{A}||_F / \sqrt{D}$ then effectively guarantees that the computed approximation $\ten{B}$ satisfies $||\ten{A}-\ten{B}||_F \leq \epsilon ||\ten{A}||_F$. 
\begin{algorithm}
\caption{Convert TT into Truncated HOSVD}
\label{alg:TT2Tucker}
\textit{\textbf{Input}}:\makebox[0pt][l]{ TT $\ten{A}^{(1)},\ldots,\ten{A}^{(D)}$ in site-1-mixed-canonical form of tensor $\ten{A}$, accuracy $\epsilon$.}\\
\textit{\textbf{Output}}:\makebox[0pt][l]{ Tucker core $\ten{S}$ in site-$D$-mixed-canonical form, orthogonal factor matrices}\\ 
\makebox[0pt][l]{$\mat{U}_1,\ldots,\mat{U}_D$ of approximation $\ten{B}$ such that $||\ten{A}-\ten{B}||_F \leq \epsilon ||\ten{A}||_F$.}
\begin{algorithmic}[1]
\STATE $\delta \gets \frac{\epsilon\,||\ten{A}||_F}{\sqrt{D}}$
\FOR{ $d=1:D$}
\STATE $\mat{U}_d,\,\mat{S},\,\mat{V}^T \gets \textrm{SVD}_{\delta}(\mat{A}_d)$ \hfill \% $\mat{A}_d = \mat{U}_d\,\mat{S}\,\mat{V}^T + \mat{E}$, $||\mat{E}||_F \leq \delta$, $S_d = \textrm{rank}(\mat{S})$.
\STATE $\ten{T} \gets \textrm{reshape}(\mat{S}\mat{V}^T,[S_d,R_d,R_{d+1}])$
\STATE $\ten{T} \gets \textrm{permute}(\ten{T},[2,1,3])$
\IF{$d==D$}
\STATE $\ten{S}^{(d)} \gets \ten{T}$
\ELSE 
%\STATE \% Construct site-$(d+1)$-mixed-canonical form 
\STATE $\mat{T} \gets \textrm{reshape}(\ten{T},[R_dS_d,R_{d+1}])$
\STATE $\mat{Q},\,\mat{R} \gets \textrm{QR}(\mat{T})$
\STATE $\ten{S}^{(d)} \gets \textrm{reshape}(\mat{Q},[R_d,S_d,R_{d+1}])$
\STATE $\ten{A}^{(d+1)} \gets \ten{A}^{(d+1)} \times_1 \mat{R}$
\ENDIF
\ENDFOR
\end{algorithmic}
\end{algorithm}

\subsection{Computational complexity}
In this subsection we briefly analyze the computational complexity of Algorithm~\ref{alg:TT2Tucker}. For notational convenience we will assume that a $D$-way tensor $\ten{A} \in \mathbb{R}^{I \times \cdots \times I}$ is represented by a TT with uniform TT-rank $R$. An additional assumption is that $I < R^2$. The computation of $D$ thin SVDs of $\mat{A}_d \in \mathbb{R}^{I \times R^2}$ in line 3 takes then $D(14R^2I^2+8I^3)$ flops~\cite[p.493]{golub13}. The QR decompositions in line 10 required for the computation of the site-$(d+1)$-mixed-canonical form require $(D-1)\left(2I^2(R^2-I/3)+4(R^4I-R^2I^2+I^3/3)\right)$ flops~\cite[p.249]{golub13} when performed with Householder transformations. In practical cases we have that $I \leq R^2$ and this implies that the total computational complexity for Algorithm~\ref{alg:TT2Tucker} is dominated by the $O(R^4I)$ term of the QR decompositions. If instead of a guaranteed relative approximation error a Tucker decomposition with given multilinear-rank is desired, then one can replace the SVD in line 3 of Algorithm~\ref{alg:TT2Tucker} by a randomized SVD~\cite{Halko2011} or an Implicitly Restarted Arnoldi Method~\cite{lehoucq1996deflation}. Also note that the actual complexity will depend heavily of the order of the indices, which is also the case with the sequentially truncated HOSVD. In practice, a heuristic that reduces the computational complexity is to permute the dimensions of the tensor $\ten{A}$ in an ascending manner prior to computing its Tucker decomposition~\cite[p.~A1041]{STHOSVD2012} as this permutation typically reduces the maximal value of $R$.

A Tucker decomposition where the Tucker core tensor is stored as a TT was first introduced in~\cite{TTTucker}. Algorithm 5~\cite[p.~611]{TTTucker} describes how such a decomposition can be obtained by means of an iterative ALS method. One disadvantage of an ALS approach, however, is that the desired TT-ranks need to be chosen a priori. An alternative DMRG approach that is able to retrieve the TT-ranks has been proposed but this comes at the cost of a computational complexity of $O(R^3I^3)$~\cite[p.~612]{TTTucker}.

\section{Tensor train to MERA}
\label{sec:TT2MERA}
\begin{figure}[tb]
    \centering
    \input{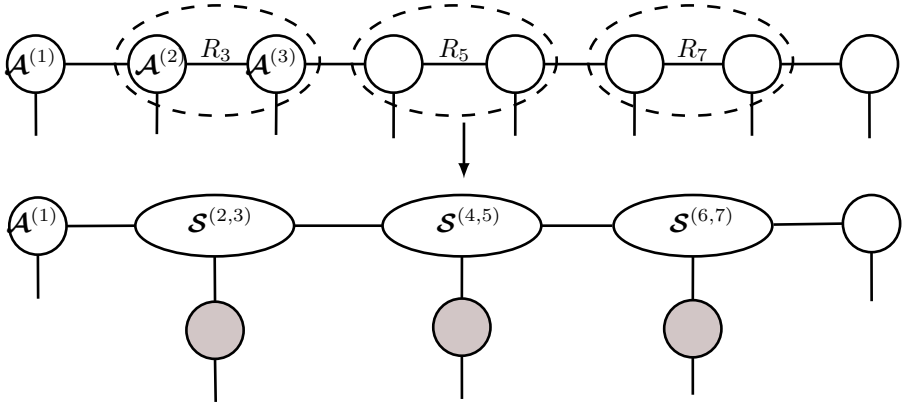}
    \caption{Diagram of disentangler computation through a HOSVD step in the TT format as described in Algorithm~\ref{alg:TT2MERA}.}
    \label{fig:disstep}
\end{figure}

\begin{figure}[tb]
    \centering
    \input{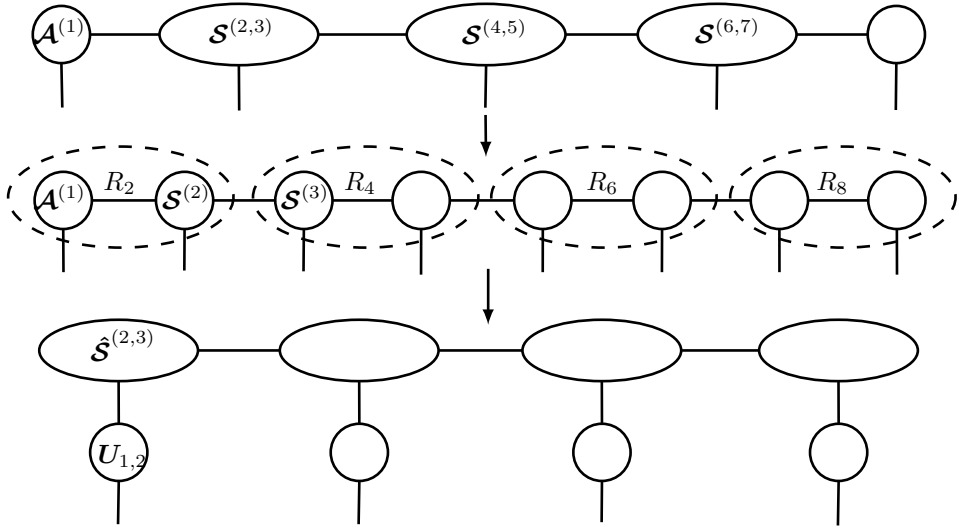}
    \caption{Diagram of isometry computation through a truncated HOSVD step in the TT format as described in Algorithm~\ref{alg:TT2MERA}.}
    \label{fig:isostep}
\end{figure}

The conversion of a TT into a MERA can be done via a sequence of HOSVD and truncated HOSVD computations. The disentanglers are computed through an HOSVD while the isometries are obtained through a truncated HOSVD. The conversion algorithm will be demonstrated through an illustrative example that consists of a TT with eight TT-cores with dimensions $R_d \times I \times R_{d+1}$ for $d=1,\ldots,8$. The goal is to compute a MERA for which the isometries convert $K=2$ indices into one. As demonstrated in Figure~\ref{fig:coarsegrainingMERA}, the ``action" of the first MERA layer is the application of three disentanglers. The diagram representation of the required operations to find these disentanglers is shown in Figure~\ref{fig:disstep}. The required disentanglers are orthogonal transformations on three index pairs. The relevant TT-cores are contracted over their auxiliary indices $R_3,R_5,R_7$ to obtain so-called ``supercores". For example, TT-cores $\ten{A}^{(2)}$ and $\ten{A}^{(3)}$ are combined into a supercore $\ten{A}^{(2,3)} \in \mathbb{R}^{R_2 \times I^2 \times R_4}$, where the two free indices of size $I$ are combined into one multi-index of size $I^2$. Algorithm~\ref{alg:TT2Tucker} is then applied to these supercores with a full SVD in order to obtain the desired disentanglers. The bottom row of Figure~\ref{fig:disstep} shows the obtained partial Tucker core with the orthogonal factor matrices, which will serve as the transposes of the disentanglers. For example, Algorithm~\ref{alg:TT2Tucker} allows us to write
\begin{align*}
    \ten{A}^{(2,3)} &= \ten{S}^{(2,3)}\; \times_2 \; \mat{U}_{2,3},
\end{align*}
where $\ten{S}^{(2,3)} \in \mathbb{R}^{R_2\times I^2 \times R_4}$ is represented by the leftmost oval of the bottom row in Figure~\ref{fig:disstep} and {$\mat{U}_{2,3} \in \mathbb{R}^{I^2 \times I^2}$} is an orthogonal matrix. The desired disentangler is then obtained by reshaping $\mat{U}_{2,3}^T$ into a cubical 4-way tensor of dimension $I$. 
The partial Tucker core is now used as the starting point for obtaining the isometries, as shown in the top row of Figure~\ref{fig:isostep}. The supercores, represented by the ovals in the top row of Figure~\ref{fig:isostep}, first need to be split back into separate TT-cores through an SVD, e.g the supercore $\ten{S}^{(2,3)}$ is reshaped into the
$R_2I \times IR_4$ matrix $\mat{S}_{2,3}$ 
\begin{align*}
    \mat{S}_{2,3}   &= \mat{U} \; \mat{S} \; \mat{V}^T,\\
                    &=\mat{S}_2 \; \mat{S}_3,
\end{align*}
with $\mat{S}_2 := \mat{U}$ and $\mat{S}_3:= \mat{S}\mat{V}^T$. The rank $R_3$ is determined as the number of nonzero singular values such that $\mat{S}_2 \in \mathbb{R}^{R_2I\times R_3}$ and $\mat{S}_3 \in \mathbb{R}^{R_3 \times IR_4}$. The desired TT-cores are obtained by reshaping $\mat{S}_2$, $\mat{S}_3$ into the desired 3-way tensors. In this way we arrive at the second row from the top of Figure~\ref{fig:isostep}. In a MERA with $K=2$ are the isometries orthogonal transformations that convert two consecutive TT indices into one index. The next step is therefore to form new supercores by summing over auxiliary indices $R_2$, $R_4$, $R_6$ and $R_8$. Applying Algorithm~\ref{alg:TT2Tucker} with a truncated SVD then results in the desired isometries. Indeed, the first supercore $\ten{\hat{A}}^{(1,2)}$ can then be written as
\begin{align*}
    \ten{\hat{A}}^{(1,2)} &= \ten{\hat{S}}^{(1,2)} \times_2 \mat{U}_{1,2},
\end{align*}
with $\ten{\hat{S}}^{(1,2)} \in \mathbb{R}^{1 \times S^2 \times R_3}$ and $\mat{U}_{1,2} \in \mathbb{R}^{I^2 \times S}$. The bottom row of Figure~\ref{fig:isostep} shows the diagram of the truncated HOSVD in TT form. The desired isometry is obtained by reshaping $\mat{U}_{1,2}$ into a $I \times I \times S$ matrix, where $S$ is the truncated index. Theorem~\ref{theo:approxerror} allows us to quantify the absolute approximation error due to truncation at each isometry step in the formation of the MERA and compute a MERA that approximates a given tensor with a guaranteed relative error. If sufficient MERA layers have been computed through this procedure, then the remaining Tucker core can be retained as the top tensor. This final step also ensures that the norm of the MERA is completely determined by the top tensor. The pseudocode for the whole algorithm is presented in Algorithm~\ref{alg:TT2MERA}.
\begin{algorithm}
\caption{MERACLE: Convert TT into a MERA}
\label{alg:TT2MERA}
\textit{\textbf{Input}}:\makebox[0pt][l]{ TT $\ten{A}^{(1)},\ldots,\ten{A}^{(D)}$ in site-$K$-mixed-canonical form, order $K$ of the isometries,}\\
\makebox[0pt][l]{accuracy $\epsilon$.}\\
\textit{\textbf{Output}}:\makebox[0pt][l]{ Isometries, disentanglers and top tensor of a MERA $\ten{B}$ such that}\\
\makebox[0pt][l]{ $||\ten{A}-\ten{B}||_F \leq \epsilon ||\ten{A}||_F$.}
\begin{algorithmic}[1]
\STATE $N \gets $ total number of isometries in the MERA.
\STATE $\delta \gets \frac{\epsilon\,||\ten{A}||_F}{\sqrt{N}}$.
\FOR{ each MERA layer}
\STATE Compute supercores of the TT according to desired disentangler locations.
\STATE Apply Algorithm~\ref{alg:TT2Tucker} with full SVD on all supercores.
\STATE Retrieve disentanglers from the orthogonal HOSVD factor matrices.
\STATE Split supercores of the partial Tucker core with the SVD.
\STATE Compute new supercores according to the location and order $K$ of the isometries.
\STATE Apply Algorithm~\ref{alg:TT2Tucker} with a truncated $\textrm{SVD}_{\delta}$ on all supercores.
\STATE Retrieve isometries from the truncated HOSVD factor matrices.
\IF {final MERA layer}
\STATE Retain single Tucker core tensor as the top tensor.
\ELSE
\STATE Split supercores of the Tucker core with an SVD.
\ENDIF
\ENDFOR
\end{algorithmic}
\end{algorithm}


\section{Iterative algorithm for finding a rank-lowering disentangler}
\label{sec:Procrustes}
Assuming that an exact low-rank MERA exists for a given TT, Algorithm~\ref{alg:TT2MERA} will typically fail to find it. In practice, the output dimensions $S$ of the isometries will simply be the product of the input dimensions $I_1I_2\cdots I_K$ and no truncation is ever performed. This leads to an exponential growth of the isometry output dimensions as a function of the number of MERA layers. The problem with Algorithm~\ref{alg:TT2MERA} is that it fails to find the correct disentanglers. In order to explain the issue at hand, we first need to explain the workings of a disentangler in a bit more detail.
\subsection{Disentangler}
\begin{figure}[tb]
    \centering
    \input{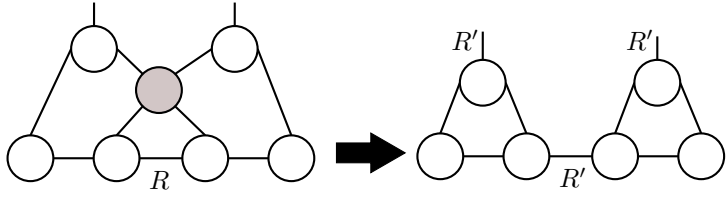}
    \caption{The disentangler reduces the TT-rank from $R$ to $R'$ with $R <R'$, allowing the two isometries to truncate to $R'$ without any loss of accuracy.}
    \label{fig:coreIdeaDisentangler}
\end{figure}

As mentioned earlier in Section~\ref{subsec:mera}, disentanglers were originally introduced in order to remove possible correlations between neighbouring indices in order to avoid high TT-ranks after coarse-graining~\cite{MERA2008}. Figure~\ref{fig:coreIdeaDisentangler} illustrates the key effect of a disentangler on a simple example of 4 TT-cores with dimensions $I_1=I_2=I_3=I_4=I$. Note that the maximal TT-rank $R$ between the second and third TT-core is $I^2$. Suppose that $R=I^2$. It is straightforward to see that having no disentangler implies that the output dimensions of the two isometries needs to be $R=I^2$, as two indices with dimension $I$ are simply combined into one multi-index. Now suppose that prior to the isometries, a disentangler can be applied to the second and third TT-cores such that the TT-rank $R$ is reduced to $R' < R$. In this case, the two isometries can truncate the dimensions $I^2$ down to $R'$ without the loss of any accuracy. Unfortunately, the disentanglers obtained from an HOSVD in Algorithm~\ref{alg:TT2MERA} do not reduce the TT-ranks, which implies that none of the isometries can effectively truncate the dimensions. If we are able to develop an algorithm that can find a rank-lowering disentangler, then Corollary~\ref{cor:multilinearrankbound} automatically guarantees that the truncated HOSVD in line 9 of Algorithm~\ref{alg:TT2MERA} will find an optimal isometry. In the next subsection we propose an iterative algorithm that attempts to recover rank-lowering disentanglers.
\subsection{Iterative orthogonal Procrustes algorithm}
Before stating the problem of finding the optimal disentangler in a formal way, we first introduce some convenient notation.
\begin{definition}
For a supercore $\ten{A}^{(d,d+1)} \in \mathbb{R}^{R_d \times I_d I_{d+1} \times R_{d+2}}$ we define the following matricizations
\begin{align*}
    \mat{A}^{(d,d+1)} &  \in \mathbb{R}^{R_dI_d \times I_{d+1}R_{d+2}}, \\
    \mat{A} & \in \mathbb{R}^{I_dI_{d+1}\times R_dR_{d+2}}.
\end{align*}
These matrices are per definition related to one another via the shuffling operator \emph{shuf} and its inverse
\begin{align*}
    \mat{A} &= \textrm{shuf}\left(\mat{A}^{(d,d+1)}\right), \\
    \mat{A}^{(d,d+1)}&= \textrm{shuf}^{-1}\left(\mat{A}\right).
\end{align*}
\end{definition}
With these definitions the optimal disentangler problem can now be formulated.
\begin{problem}
\label{problem:procrustes}
Given a supercore $\ten{A}^{(d,d+1)}$ for which \mbox{rank$\left(\mat{A}^{(d,d+1)}\right)=R$}, find an orthogonal matrix $\mat{V} \in \mathbb{R}^{I_d I_{d+1} \times I_d I_{d+1}}$ such that
\begin{align*}
    \mat{A}' &:= \mat{V}\, \textrm{shuf}\left(\mat{A}^{(d,d+1)}\right) =\mat{V}\, \mat{A}
\end{align*}
with rank$\left(\textrm{shuf}^{-1}\left(\mat{A}'\right)\right)=R' < R$.
\end{problem}
Problem~\ref{problem:procrustes} is essentially an orthogonal Procrustes problem in $\mat{A}$ with the additional constraint that the orthogonal transformation $\mat{V}$ lowers the rank of $\mat{A}^{(d,d+1)}$. The difficulty is that both $\mat{A}'$ and $\textrm{shuf}^{-1}\left(\mat{A}'\right)$ are unknown. We therefore propose to solve the orthogonal Procrustes problem in an iterative manner, where we fix $\mat{A}^{(k,k+1)'}$ in every iteration to a low-rank approximation of $\mat{A}^{(d,d+1)}$. The computational complexity of solving the orthogonal Procrustes problem every iteration is $O((I_dI_{d+1})^3)$, as this amounts to computing the SVD of $\mat{A}'\,\mat{A}^T$. The proposed iterative algorithm is presented in pseudocode as Algorithm~\ref{alg:disentangler}. The stopping criterion can be set to a fixed maximum number of iterations or one can inspect the rank-gap $\sigma_{R'} / \sigma_{R'+1}$ of $\mat{A}^{(d,d+1)}$ and stop the iterations as soon as this gap has reached a certain order of magnitude. A low-rank approximation of $\mat{A}^{(d,d+1)}$ can be computed via its SVD. At this moment, there is no formal proof of convergence for Algorithm~\ref{alg:disentangler}, nor is it known what the conditions for convergence are. The best that we are currently able to do is to empirically show the successful application of this algorithm and to explore its properties based on extensive numerical experiments.
\begin{algorithm}
\caption{Iterative disentangler computation}
\label{alg:disentangler}
\textit{\textbf{Input}}:\makebox[0pt][l]{ supercore $\ten{A}^{(d,d+1)}  \in \mathbb{R}^{R_d \times I_d I_{d+1}\times R_{d+2}}$.}\\
\textit{\textbf{Output}}:\makebox[0pt][l]{ Disentangler $\ten{V}$ that reduces the TT-rank to $R'$.}
\begin{algorithmic}
\STATE 1. $\mat{V} \gets \mat{I}$ \hfill \% initialize with identity matrix
\WHILE{stopping criterion not true}
\STATE 2. $\mat{A}^{(d,d+1)'} \gets $ low-rank approximation of $\mat{A}^{(d,d+1)}$.
\STATE 3. $\mat{A} \gets \textrm{shuf}\left(\mat{A}^{(d,d+1)} \right)$ and $\mat{A}' \gets \textrm{shuf}\left(\mat{A}^{(d,d+1)'} \right)$.
\STATE 4. $\hat{\mat{V}} \gets $ solve orthogonal Procrustes problem that minimizes $||\mat{V}\mat{A}-\mat{A}'||_F^2$.
\STATE 5. $\mat{A} \gets \hat{\mat{V}}\,\mat{A}$.
\STATE 6. $\mat{A}^{(d,d+1)} \gets \textrm{shuf}^{-1}\left(\mat{A} \right)$.
\STATE 7. $\mat{V} \gets \hat{\mat{V}}\,\mat{V}$.
\ENDWHILE
\STATE 8. $\ten{V} \gets \textrm{reshape}(\mat{V},[I_d,I_{d+1},I_d,I_{d+1}])$.
\end{algorithmic}
\end{algorithm}

%\section{Truncation of MERA-ranks}
%\begin{algorithm}
%\caption{Truncation of MERA-ranks}
%\label{alg:Mera_round}
%\textit{\textbf{Input}}:\makebox[0pt][l]{ MERA of a tensor $\ten{A}$, prescribed %accuracy $\epsilon$.}\\
%\textit{\textbf{Output}}:\makebox[0pt][l]{ MERA-approximation $\ten{B}$ of $\ten{A}$ with reduced ranks such that \\
%$||\ten{A}-\ten{B}||_F \leq \epsilon\, ||\ten{A}||_F$.}
%\begin{algorithmic}
%\STATE 1. $\ten{T}^{(1)},\ldots,\ten{T}^{(d)} \gets$ Convert top MERA tensor %$\ten{T}$ in TT form through the TT-SVD algorithm.
%\STATE 2. Apply the TT-rounding algorithm on the TT of $\ten{T}$ with prescribed %accuracy $\epsilon$.
%\STATE 3. Compute a THOSVD of the TT with Algorithm~\ref{alg:TT2Tucker}.
%\STATE 4. Propagate the orthogonal factor matrices of the THOSVD through the MERA %structure.
%\end{algorithmic}
%\end{algorithm}

\section{Experiments}
\label{sec:experiments}
In this section we demonstrate the computational efficiency of Algorithms~\ref{alg:TT2Tucker},~\ref{alg:TT2MERA} and~\ref{alg:disentangler} through numerical experiments. All algorithms were implemented in MATLAB and the experiments were performed on a desktop computer with a 4-core processor running at 3.6~GHz with 16~GB RAM.

\subsection{Converting a TT into Tucker - compression of simulation results}
In this experiment we demonstrate Algorithm~\ref{alg:TT2Tucker} and how a representation of a Tucker decomposition can benefit compression without loss of accuracy. Inspired by the example discussed in~\cite[p.~A1047]{STHOSVD2012}, a tensor decomposition is used for the compression the solution $u(x,y,t)$ of
\begin{align*}
    \frac{\partial u}{\partial t} &=\frac{\partial^2 u}{\partial x^2} + \frac{\partial^2 u}{\partial y^2} 
\end{align*}
on the unit square $[0,1]^2$ with boundary condition $0.25-|0.5-x|\cdot |0.5-y|$, which also describes the initial temperature distribution over the entire square. The PDE was discretized with a uniform mesh with cell size $(\Delta s,\Delta s, \Delta t)$ and solved with the explicit Euler method, using a time step $0.25\Delta s^2$ to ensure numerical stability. We set $\Delta s=10^{-2}$ and $\Delta t=0.25\cdot 10^{-4}$ and simulate for about $0.25$ seconds, resulting in a tensor of size $100\times 100 \times 10000$. The upper bound on the relative approximation error when computing tensor decompositions is set to $10^{-3}$. We compare the sequentially truncated HOSVD (STHOSVD) with both the TT and Tucker decomposition in TT form. The STHOSVD is computed with the mlsvd command of the Tensorlab toolbox~\cite{tensorlab3.0}, while the conversion of the original data tensor into a TT is done via the TT-SVD algorithm~\cite[p.~2301]{ivanTT}. The TT is converted into a Tucker decomposition via Algorithm~\ref{alg:TT2Tucker}. We consider two cases. In the first case, we compute the three tensor decompositions on the original solution tensor, while in the second case we first reshape the original data into a 16-way tensor by factorization of all dimensions into their prime components. All results are shown in Table~\ref{tab:STHOSVDvsTT}. The compression column contains the ratio between how many numbers are required to store the original tensor and how many numbers are required to store the decomposition. Not much difference in neither the total runtime, relative error or compression can be observed when the simulation solution is kept as a 3-way tensor. The computation of an STHOSVD of the 16-way tensor takes about 3 times longer than computing the corresponding TT. The resulting decomposition is also not able to compress the data very much as each of the dimensions of the 16-way tensor consist of (small) prime factors. The TT and Tucker decomposition in TT form, however, result in a saving of around 12,000, which is an improvement of more than 10 times compared to the 3-way case. The time required for Algorithm~\ref{alg:TT2Tucker} to compute the Tucker decomposition was in both cases negligible compared to the runtime of the TT-SVD algorithm.

\begin{table}[ht]
\centering
\caption{Comparison of different tensor decompositions in compressing the results from a numerical simulation of the 2D heat equation.}%\vspace{3pt}
\label{tab:STHOSVDvsTT} % is used to refer this table in the text
\begin{tabular}{@{}lrrr@{}}
                    & Time (s)  & Relative error    & Compression \\\midrule
STHOSVD (3-way)    & $3.279$     & $4.74\cdot 10^{-4}$ & $1055$     \\
TT (3-way)         & $3.579$     & $8.54\cdot 10^{-4}$ & $1101$     \\
TT-Tucker (3-way)   & $0.019$     & $8.54\cdot 10^{-4}$ & $1116$\\
STHOSVD (16-way)    & $26.71$     & $6.68\cdot 10^{-4}$ & $3.255$     \\
TT (16-way)         & $8.325$     & $6.11\cdot 10^{-4}$ & $12,572$     \\
TT-Tucker (16-way)   & $0.005$     & $6.11\cdot 10^{-4}$ & $12,229$
\end{tabular}
\end{table}

\subsection{Comparison of HOSVD with Algorithm~\ref{alg:disentangler}}
\label{subsec:exp2}
\begin{figure}[tb]
    \centering
    % Graphic for TeX using PGF
% Title: /home/kbatselier/Dropbox/Figures/Dia/Image_disentanglement_example.dia
% Creator: Dia v0.97+git
% CreationDate: Thu Nov 21 15:54:50 2019
% For: kbatselier
% \usepackage{tikz}
% The following commands are not supported in PSTricks at present
% We define them conditionally, so when they are implemented,
% this pgf file will use them.
\ifx\du\undefined
  \newlength{\du}
\fi
\setlength{\du}{4\unitlength}
\begin{tikzpicture}[even odd rule]
\pgftransformxscale{1.000000}
\pgftransformyscale{-1.000000}
\definecolor{dialinecolor}{rgb}{0.000000, 0.000000, 0.000000}
\pgfsetstrokecolor{dialinecolor}
\pgfsetstrokeopacity{1.000000}
\definecolor{diafillcolor}{rgb}{1.000000, 1.000000, 1.000000}
\pgfsetfillcolor{diafillcolor}
\pgfsetfillopacity{1.000000}
\pgfsetlinewidth{0.200000\du}
\pgfsetdash{}{0pt}
\pgfsetbuttcap
\pgfsetmiterjoin
\pgfsetlinewidth{0.200000\du}
\pgfsetbuttcap
\pgfsetmiterjoin
\pgfsetdash{}{0pt}
\definecolor{diafillcolor}{rgb}{1.000000, 1.000000, 1.000000}
\pgfsetfillcolor{diafillcolor}
\pgfsetfillopacity{1.000000}
\pgfpathellipse{\pgfpoint{40.850000\du}{10.325000\du}}{\pgfpoint{3.475000\du}{0\du}}{\pgfpoint{0\du}{3.475000\du}}
\pgfusepath{fill}
\definecolor{dialinecolor}{rgb}{0.000000, 0.000000, 0.000000}
\pgfsetstrokecolor{dialinecolor}
\pgfsetstrokeopacity{1.000000}
\pgfpathellipse{\pgfpoint{40.850000\du}{10.325000\du}}{\pgfpoint{3.475000\du}{0\du}}{\pgfpoint{0\du}{3.475000\du}}
\pgfusepath{stroke}
\pgfsetlinewidth{0.200000\du}
\pgfsetdash{}{0pt}
\pgfsetbuttcap
\pgfsetmiterjoin
\pgfsetlinewidth{0.200000\du}
\pgfsetbuttcap
\pgfsetmiterjoin
\pgfsetdash{}{0pt}
\definecolor{diafillcolor}{rgb}{1.000000, 1.000000, 1.000000}
\pgfsetfillcolor{diafillcolor}
\pgfsetfillopacity{1.000000}
\pgfpathellipse{\pgfpoint{56.825000\du}{19.285000\du}}{\pgfpoint{3.475000\du}{0\du}}{\pgfpoint{0\du}{3.475000\du}}
\pgfusepath{fill}
\definecolor{dialinecolor}{rgb}{0.000000, 0.000000, 0.000000}
\pgfsetstrokecolor{dialinecolor}
\pgfsetstrokeopacity{1.000000}
\pgfpathellipse{\pgfpoint{56.825000\du}{19.285000\du}}{\pgfpoint{3.475000\du}{0\du}}{\pgfpoint{0\du}{3.475000\du}}
\pgfusepath{stroke}
\pgfsetlinewidth{0.200000\du}
\pgfsetdash{}{0pt}
\pgfsetbuttcap
{
\definecolor{diafillcolor}{rgb}{0.000000, 0.000000, 0.000000}
\pgfsetfillcolor{diafillcolor}
\pgfsetfillopacity{1.000000}
% was here!!!
\definecolor{dialinecolor}{rgb}{0.000000, 0.000000, 0.000000}
\pgfsetstrokecolor{dialinecolor}
\pgfsetstrokeopacity{1.000000}
\draw (44.325000\du,10.325000\du)--(56.825000\du,15.810000\du);
}
\pgfsetlinewidth{0.200000\du}
\pgfsetdash{}{0pt}
\pgfsetbuttcap
\pgfsetmiterjoin
\pgfsetlinewidth{0.200000\du}
\pgfsetbuttcap
\pgfsetmiterjoin
\pgfsetdash{}{0pt}
\definecolor{diafillcolor}{rgb}{1.000000, 1.000000, 1.000000}
\pgfsetfillcolor{diafillcolor}
\pgfsetfillopacity{1.000000}
\pgfpathellipse{\pgfpoint{40.850000\du}{29.080000\du}}{\pgfpoint{3.475000\du}{0\du}}{\pgfpoint{0\du}{3.475000\du}}
\pgfusepath{fill}
\definecolor{dialinecolor}{rgb}{0.000000, 0.000000, 0.000000}
\pgfsetstrokecolor{dialinecolor}
\pgfsetstrokeopacity{1.000000}
\pgfpathellipse{\pgfpoint{40.850000\du}{29.080000\du}}{\pgfpoint{3.475000\du}{0\du}}{\pgfpoint{0\du}{3.475000\du}}
\pgfusepath{stroke}
\pgfsetlinewidth{0.200000\du}
\pgfsetdash{}{0pt}
\pgfsetbuttcap
{
\definecolor{diafillcolor}{rgb}{0.000000, 0.000000, 0.000000}
\pgfsetfillcolor{diafillcolor}
\pgfsetfillopacity{1.000000}
% was here!!!
\definecolor{dialinecolor}{rgb}{0.000000, 0.000000, 0.000000}
\pgfsetstrokecolor{dialinecolor}
\pgfsetstrokeopacity{1.000000}
\draw (44.325000\du,29.080000\du)--(56.825000\du,22.760000\du);
}
\pgfsetlinewidth{0.200000\du}
\pgfsetdash{}{0pt}
\pgfsetbuttcap
\pgfsetmiterjoin
\pgfsetlinewidth{0.200000\du}
\pgfsetbuttcap
\pgfsetmiterjoin
\pgfsetdash{}{0pt}
\definecolor{diafillcolor}{rgb}{0.749020, 0.749020, 0.749020}
\pgfsetfillcolor{diafillcolor}
\pgfsetfillopacity{1.000000}
\pgfpathellipse{\pgfpoint{28.650000\du}{19.680000\du}}{\pgfpoint{3.475000\du}{0\du}}{\pgfpoint{0\du}{3.475000\du}}
\pgfusepath{fill}
\definecolor{dialinecolor}{rgb}{0.000000, 0.000000, 0.000000}
\pgfsetstrokecolor{dialinecolor}
\pgfsetstrokeopacity{1.000000}
\pgfpathellipse{\pgfpoint{28.650000\du}{19.680000\du}}{\pgfpoint{3.475000\du}{0\du}}{\pgfpoint{0\du}{3.475000\du}}
\pgfusepath{stroke}
\pgfsetlinewidth{0.200000\du}
\pgfsetdash{}{0pt}
\pgfsetbuttcap
{
\definecolor{diafillcolor}{rgb}{0.000000, 0.000000, 0.000000}
\pgfsetfillcolor{diafillcolor}
\pgfsetfillopacity{1.000000}
% was here!!!
\definecolor{dialinecolor}{rgb}{0.000000, 0.000000, 0.000000}
\pgfsetstrokecolor{dialinecolor}
\pgfsetstrokeopacity{1.000000}
\draw (17.200000\du,6.900000\du)--(40.850000\du,6.850000\du);
}
\pgfsetlinewidth{0.200000\du}
\pgfsetdash{}{0pt}
\pgfsetbuttcap
{
\definecolor{diafillcolor}{rgb}{0.000000, 0.000000, 0.000000}
\pgfsetfillcolor{diafillcolor}
\pgfsetfillopacity{1.000000}
% was here!!!
\definecolor{dialinecolor}{rgb}{0.000000, 0.000000, 0.000000}
\pgfsetstrokecolor{dialinecolor}
\pgfsetstrokeopacity{1.000000}
\draw (31.550492\du,17.751983\du)--(40.850000\du,13.800000\du);
}
\pgfsetlinewidth{0.200000\du}
\pgfsetdash{}{0pt}
\pgfsetbuttcap
{
\definecolor{diafillcolor}{rgb}{0.000000, 0.000000, 0.000000}
\pgfsetfillcolor{diafillcolor}
\pgfsetfillopacity{1.000000}
% was here!!!
\definecolor{dialinecolor}{rgb}{0.000000, 0.000000, 0.000000}
\pgfsetstrokecolor{dialinecolor}
\pgfsetstrokeopacity{1.000000}
\draw (31.357248\du,21.869573\du)--(40.850000\du,25.605000\du);
}
\pgfsetlinewidth{0.200000\du}
\pgfsetdash{}{0pt}
\pgfsetbuttcap
{
\definecolor{diafillcolor}{rgb}{0.000000, 0.000000, 0.000000}
\pgfsetfillcolor{diafillcolor}
\pgfsetfillopacity{1.000000}
% was here!!!
\definecolor{dialinecolor}{rgb}{0.000000, 0.000000, 0.000000}
\pgfsetstrokecolor{dialinecolor}
\pgfsetstrokeopacity{1.000000}
\draw (17.297333\du,17.619552\du)--(25.834312\du,17.616162\du);
}
\pgfsetlinewidth{0.200000\du}
\pgfsetdash{}{0pt}
\pgfsetbuttcap
{
\definecolor{diafillcolor}{rgb}{0.000000, 0.000000, 0.000000}
\pgfsetfillcolor{diafillcolor}
\pgfsetfillopacity{1.000000}
% was here!!!
\definecolor{dialinecolor}{rgb}{0.000000, 0.000000, 0.000000}
\pgfsetstrokecolor{dialinecolor}
\pgfsetstrokeopacity{1.000000}
\draw (17.333242\du,21.820851\du)--(25.901744\du,21.801591\du);
}
\pgfsetlinewidth{0.200000\du}
\pgfsetdash{}{0pt}
\pgfsetbuttcap
{
\definecolor{diafillcolor}{rgb}{0.000000, 0.000000, 0.000000}
\pgfsetfillcolor{diafillcolor}
\pgfsetfillopacity{1.000000}
% was here!!!
\definecolor{dialinecolor}{rgb}{0.000000, 0.000000, 0.000000}
\pgfsetstrokecolor{dialinecolor}
\pgfsetstrokeopacity{1.000000}
\draw (17.050000\du,32.600000\du)--(40.850000\du,32.555000\du);
}
% setfont left to latex
\definecolor{dialinecolor}{rgb}{0.000000, 0.000000, 0.000000}
\pgfsetstrokecolor{dialinecolor}
\pgfsetstrokeopacity{1.000000}
\definecolor{diafillcolor}{rgb}{0.000000, 0.000000, 0.000000}
\pgfsetfillcolor{diafillcolor}
\pgfsetfillopacity{1.000000}
\node[anchor=base west,inner sep=0pt,outer sep=0pt,color=dialinecolor] at (48.325000\du,11.600000\du){$R'$};
% setfont left to latex
\definecolor{dialinecolor}{rgb}{0.000000, 0.000000, 0.000000}
\pgfsetstrokecolor{dialinecolor}
\pgfsetstrokeopacity{1.000000}
\definecolor{diafillcolor}{rgb}{0.000000, 0.000000, 0.000000}
\pgfsetfillcolor{diafillcolor}
\pgfsetfillopacity{1.000000}
\node[anchor=base west,inner sep=0pt,outer sep=0pt,color=dialinecolor] at (48.325000\du,24.525000\du){R'};
% setfont left to latex
\definecolor{dialinecolor}{rgb}{0.000000, 0.000000, 0.000000}
\pgfsetstrokecolor{dialinecolor}
\pgfsetstrokeopacity{1.000000}
\definecolor{diafillcolor}{rgb}{0.000000, 0.000000, 0.000000}
\pgfsetfillcolor{diafillcolor}
\pgfsetfillopacity{1.000000}
\node[anchor=base west,inner sep=0pt,outer sep=0pt,color=dialinecolor] at (17.829650\du,5.930000\du){$I$};
% setfont left to latex
\definecolor{dialinecolor}{rgb}{0.000000, 0.000000, 0.000000}
\pgfsetstrokecolor{dialinecolor}
\pgfsetstrokeopacity{1.000000}
\definecolor{diafillcolor}{rgb}{0.000000, 0.000000, 0.000000}
\pgfsetfillcolor{diafillcolor}
\pgfsetfillopacity{1.000000}
\node[anchor=base west,inner sep=0pt,outer sep=0pt,color=dialinecolor] at (17.829650\du,15.175000\du){I};
% setfont left to latex
\definecolor{dialinecolor}{rgb}{0.000000, 0.000000, 0.000000}
\pgfsetstrokecolor{dialinecolor}
\pgfsetstrokeopacity{1.000000}
\definecolor{diafillcolor}{rgb}{0.000000, 0.000000, 0.000000}
\pgfsetfillcolor{diafillcolor}
\pgfsetfillopacity{1.000000}
\node[anchor=base west,inner sep=0pt,outer sep=0pt,color=dialinecolor] at (33.153183\du,16.078152\du){I};
% setfont left to latex
\definecolor{dialinecolor}{rgb}{0.000000, 0.000000, 0.000000}
\pgfsetstrokecolor{dialinecolor}
\pgfsetstrokeopacity{1.000000}
\definecolor{diafillcolor}{rgb}{0.000000, 0.000000, 0.000000}
\pgfsetfillcolor{diafillcolor}
\pgfsetfillopacity{1.000000}
\node[anchor=base west,inner sep=0pt,outer sep=0pt,color=dialinecolor] at (33.153183\du,25.363869\du){I};
% setfont left to latex
\definecolor{dialinecolor}{rgb}{0.000000, 0.000000, 0.000000}
\pgfsetstrokecolor{dialinecolor}
\pgfsetstrokeopacity{1.000000}
\definecolor{diafillcolor}{rgb}{0.000000, 0.000000, 0.000000}
\pgfsetfillcolor{diafillcolor}
\pgfsetfillopacity{1.000000}
\node[anchor=base west,inner sep=0pt,outer sep=0pt,color=dialinecolor] at (17.829650\du,24.800802\du){I};
% setfont left to latex
\definecolor{dialinecolor}{rgb}{0.000000, 0.000000, 0.000000}
\pgfsetstrokecolor{dialinecolor}
\pgfsetstrokeopacity{1.000000}
\definecolor{diafillcolor}{rgb}{0.000000, 0.000000, 0.000000}
\pgfsetfillcolor{diafillcolor}
\pgfsetfillopacity{1.000000}
\node[anchor=base west,inner sep=0pt,outer sep=0pt,color=dialinecolor] at (17.829650\du,31.430000\du){I};
\end{tikzpicture}
    \caption{The disentangler reduces the TT-rank from $R$ to $R'$ with $R > R'$, allowing the two isometries to truncate to $R'$ without any loss of accuracy.}
    \label{fig:image}
\end{figure}
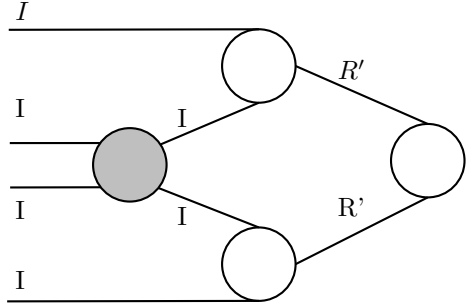

\begin{figure}[tb]
\centering
\subfigure[Top tensor image.]{
\begin{minipage}[t]{0.31\linewidth}
\includegraphics[width=\textwidth]{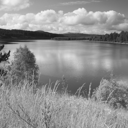}
\end{minipage}}
\subfigure[Image after isometries.]{
\begin{minipage}[t]{0.31\linewidth}
\includegraphics[width=\textwidth]{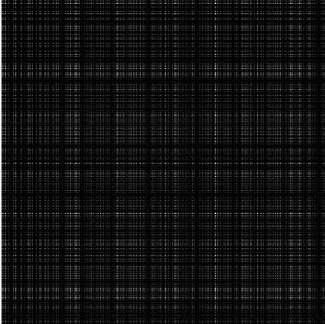}
\end{minipage}}
\subfigure[Image after disentangler.]{
\begin{minipage}[t]{0.31\linewidth}
\includegraphics[width=\textwidth]{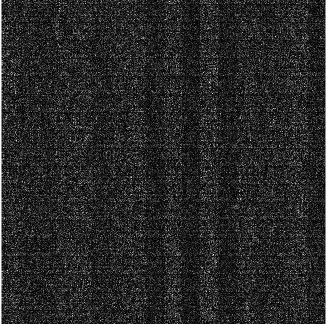}
\end{minipage}}
\caption{Image after consecutive application of isometries and a disentangler.}
\label{fig:MERAimage}
\end{figure}

\begin{figure}[tb]
\centering
\subfigure[Convergence of singular values during Algorithm 5.1.]{
\begin{minipage}[b]{0.47\linewidth}
\includegraphics[width=\textwidth]{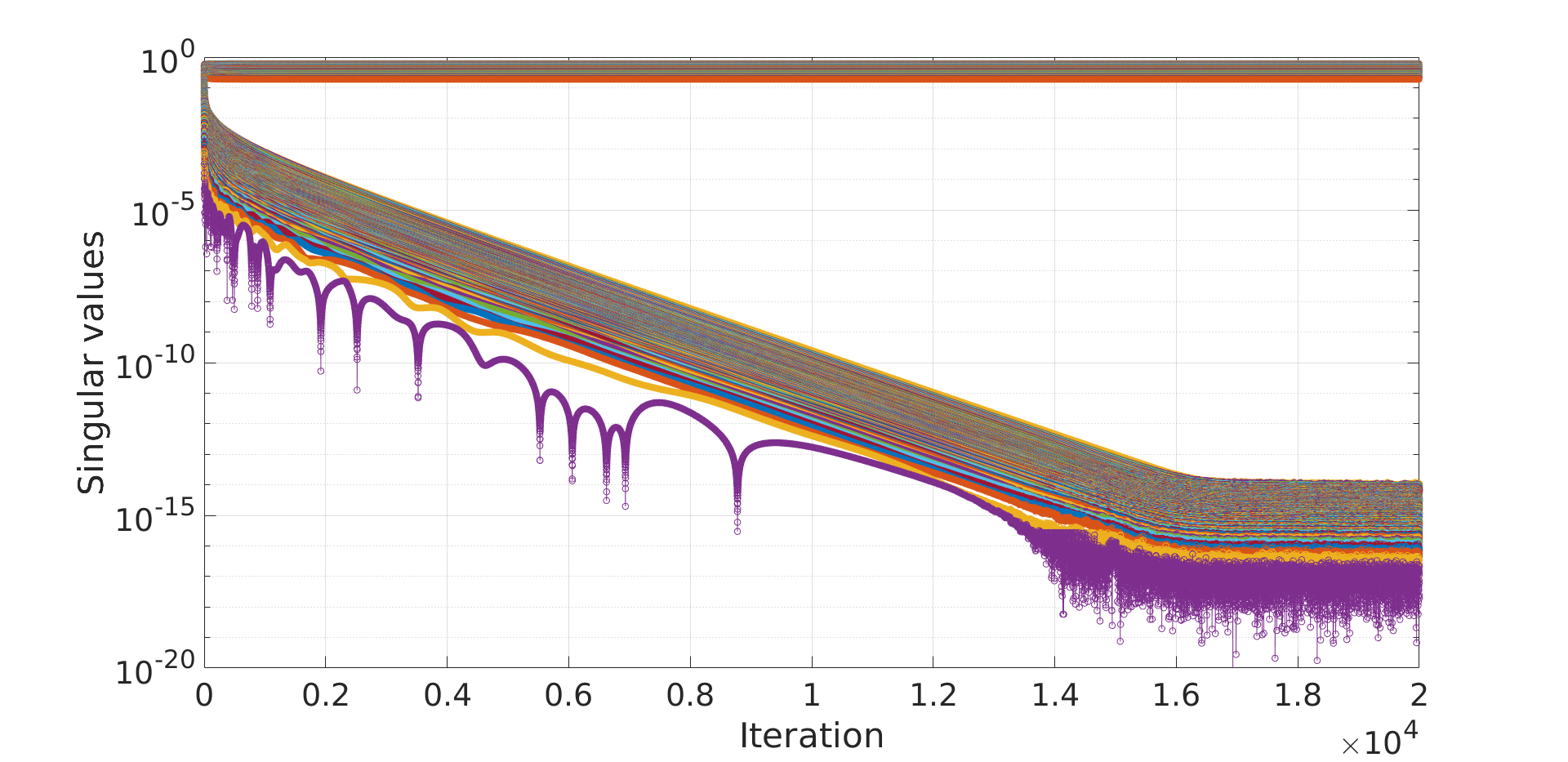}
\end{minipage}}
\subfigure[Singular value decay of $\mat{C}^{(2,3)}$ after application of disentanglers.]{
\begin{minipage}[b]{0.47\linewidth}
\includegraphics[width=\textwidth]{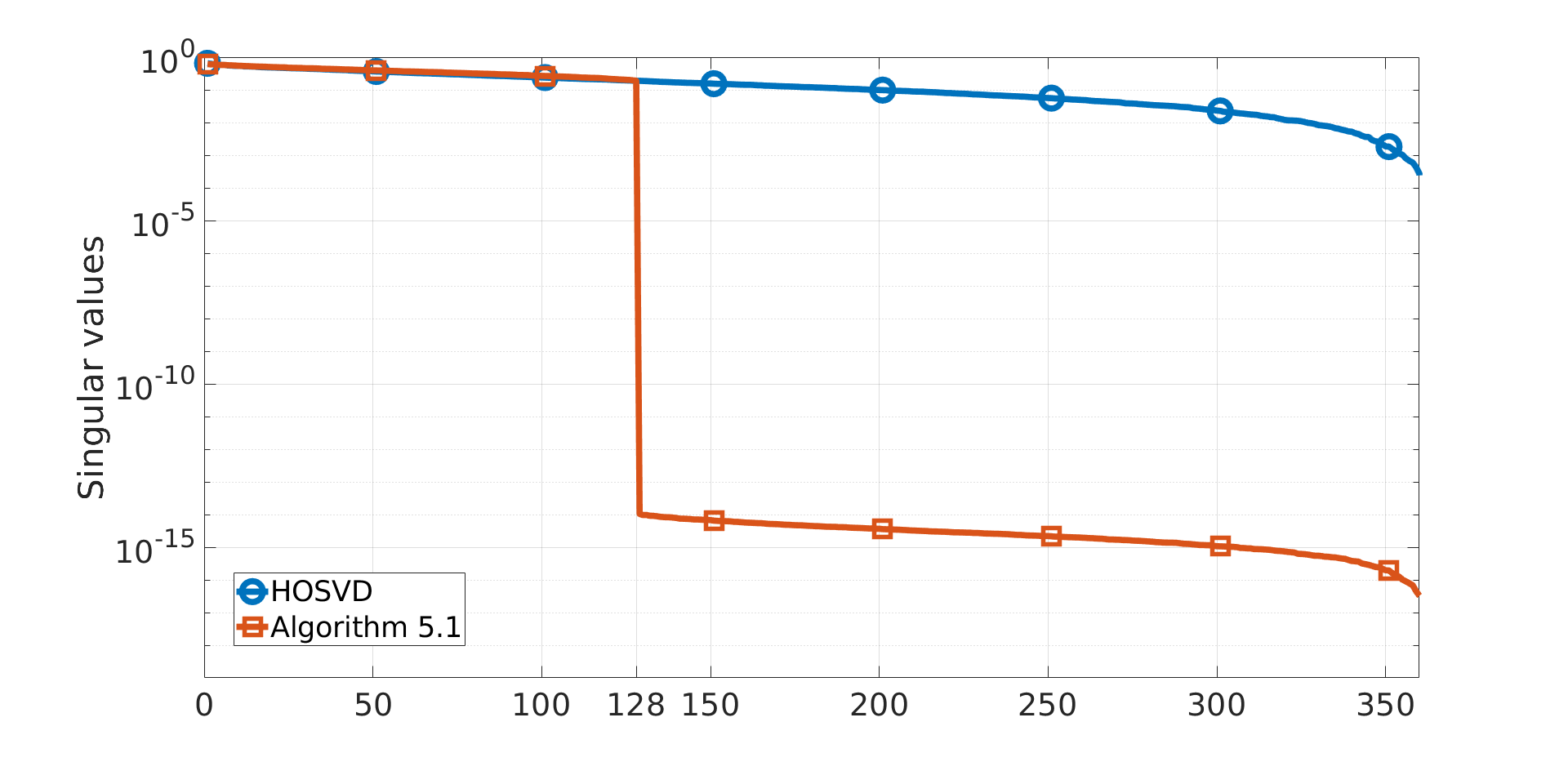}
\end{minipage}}
\caption{Singular value graphs.}
\label{fig:HOSVDvsAlgo5}
\end{figure}

In this experiment we compare Algorithm~\ref{alg:TT2Tucker} with Algorithm~\ref{alg:disentangler} to retrieve a rank-lowering disentangler. For this we consider the MERA consisting of a single layer as depicted in Figure~\ref{fig:image}. The top tensor is taken to be an $R' \times R'$ grayscale image\footnote{The image was taken from \url{http://absfreepic.com/free-photos/download/landscape-with-lake-4412x2941_12692.html}, cropped and scaled to appropriate dimensions.}. In this particular case we set $R'=128$. Both the  $I \times I \times R'$ isometry tensors $\ten{W}$ and $I \times I \times I \times I$ disentangler tensor $\ten{V}$ are found from the orthogonalization of random matrices with appropriate sizes. For this experiment we set $I=19$ and choose the isometries to be identical. Given the $R' \times R'$ grayscale image top tensor $\mat{A}$ shown in Figure~\ref{fig:MERAimage}(a), we can now apply the MERA `backwards'. The application of the two isometries on $\mat{A}$ is
\begin{align*}
    \mat{B} &= \mat{W}\,\mat{A}\,\mat{W}^T,
\end{align*}
resulting in an $I^2 \times I^2$ image $\mat{B}$, shown in Figure~\ref{fig:MERAimage}(b). The corresponding TT of the image $\mat{B}$ has TT-ranks $19$ and $128=R'$. The application of the disentangler is then performed from the following steps
\begin{align*}
& \ten{B} := \textrm{reshape}(\mat{B},[I,I,I,I]),\\    
& \ten{B}_p:= \textrm{permute}(\ten{B},[2,3,4,1]),\\
& \mat{\tilde{B}}:= \textrm{reshape}(\ten{B}_p,[I^2,I^2]),\\
& \mat{C}_p:= \mat{V}\,\mat{\tilde{B}},\\
& \ten{C}_p:= \textrm{reshape}(\mat{C}_p,[I,I,I,I]),\\
& \ten{C}:= \textrm{permute}(\ten{C}_p,[4,1,2,3]),\\
& \mat{C} := \textrm{reshape}(\ten{C},[I^2,I^2]),
\end{align*}
resulting in the $I^2 \times I^2$ image $\mat{C}$ shown in Figure~\ref{fig:MERAimage}(c). The corresponding TT of the image $\mat{C}$ has TT-ranks $19$ and $361=I^2$. This increase of the TT-rank is reflected in the image as being much more 'noisy' while the low-rank image of Figure~\ref{fig:MERAimage}(b) has a particular block structure pattern. We now compare the use of Algorithm~\ref{alg:TT2Tucker} with Algorithm~\ref{alg:disentangler} for retrieving a disentangler that is able to reduce the maximal TT-rank from 361 down to 128. Algorithm~\ref{alg:disentangler} is run on the TT of $\ten{C}$ in site-4-mixed-canonical form and a rank-128 approximation of $\mat{C}^{(2,3)}$ is used. Each iteration of Algorithm~\ref{alg:disentangler} took 0.03 seconds and, as shown in Figure~\ref{fig:HOSVDvsAlgo5}(a), about 16,000 iterations were required for the 233 smallest singular values to converge to values of about $10^{-15}$. The computed disentanglers are then applied to the supercore $\ten{C}^{(2,3)}$. The singular value decay of each corresponding $\mat{C}^{(2,3)}$ is shown in Figure~\ref{fig:HOSVDvsAlgo5}(b), where it can be clearly seen that Algorithm~\ref{alg:disentangler} is able to retrieve a disentangler that lowers the rank to the minimal value of 128.

\subsection{Limitations of Algorithm~\ref{alg:disentangler}}
We revisit the example from subsection~\ref{subsec:exp2} and explore the validity of Algorithm~\ref{alg:disentangler} for different values of $R'$ and $I$, as it is yet unclear under which conditions we are able to retrieve an exact rank-lowering disentangler. If $I$ is fixed, then the rank of $\mat{C}$ is $I^2$ for the particular MERA of Section~\ref{subsec:exp2} and it appears that there exists a minimal value $R'_{\textrm{min}}$ such that Algorithm~\ref{alg:disentangler} does not converge for values $R' < R'_{\textrm{min}}$. There is however an exception to this observation in that Algorithm~\ref{alg:disentangler} always converges if $R'=1$. Table~\ref{tab:Algo5fail} lists all values of $R'_{\textrm{min}}$ for values of $I$ going from 2 up to 14, where convergence of Algorithm~\ref{alg:disentangler} was determined from inspecting the singular value decay as in Figure~\ref{fig:HOSVDvsAlgo5}(a). A first observation is that $R'_{\textrm{min}}$ grows slowly compared to $R=I^2$, which implies that the range of values of $R'$ for which Algorithm~\ref{alg:disentangler} converges gets larger as $I$ grows. The reason for the existence of this $R'_{\textrm{min}}$ is yet to be fully understood.\\

\begin{table}[t]
    \centering
    \caption{Minimal value of  $R'$ for which Algorithm~\ref{alg:disentangler} converges as a function of $I$.}
    \label{tab:Algo5fail}    
    \begin{tabular}{l|r r r r r r r r r r r r r}
    $I$  & 2 & 3 & 4 & 5 & 6 & 7 & 8 & 9 & 10 & 11 & 12 & 13& 14   \\  \hline
$R'_{\textrm{min}}$  & 2 & 4 & 6 & 9 &12 & 16  & 20  & 25 & 30 & 37  & 44 & 51 & 59    
    \end{tabular}
\end{table}

\begin{figure}
    \centering
    \includegraphics[width=.8\textwidth]{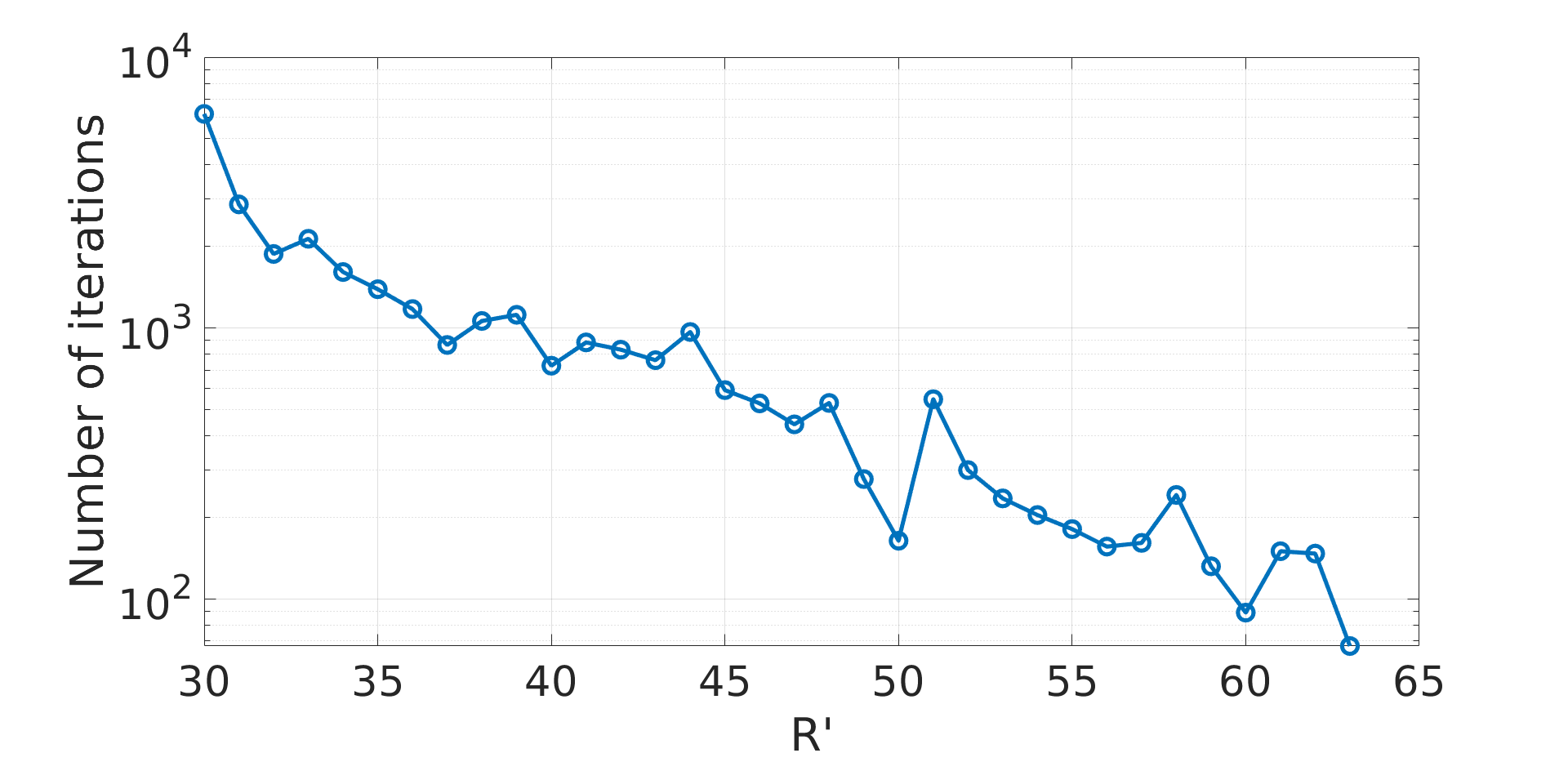}
    \caption{The number of iterations required for Algorithm~\ref{alg:disentangler} to reach a rank-gap of $\sigma_{R'}/\sigma_{R'+1} = 10^{12}$ as a function of $R'$ when $I=8$.}
    \label{fig:IterationsvsR}
\end{figure}

A second observation relates to the rate of convergence. It turns out that Algorithm~\ref{alg:disentangler} converges faster as the difference between $R$ and $R'$ becomes smaller. This is illustrated in Figure~\ref{fig:IterationsvsR} where the number of iterations required for Algorithm~\ref{alg:disentangler} to reach a rank-gap of $\sigma_{R'}/\sigma_{R'+1} = 10^{12}$ is shown for varying $R'$ when $I=8$. An approximately exponential growth in the number of required iterations can be seen as the difference of $R'$ with $R=8^2=64$ grows larger. This exponential growth might explain the existence of $R'_{\textrm{min}}$ as a value of $R'$ for which convergence becomes `infinitely slow'. These observations will serve as a starting point to investigate the exact nature of why and when Algorithm~\ref{alg:disentangler} works, apart from the empirical study herein.

\subsection{Comparison of compression capability between a TT and a MERA on a large-scale example}
\label{sec:exp3}
In this experiment we compare the compression capability between a TT and a MERA. We also apply Algorithm~\ref{alg:TT2MERA} on a large-scale example for which a 12-way cubical tensor $\ten{A}$ of dimension 10 is generated that is exactly represented by a 2-layer MERA, where each of the isometries reduces $K=2$ indices into 1 index $S=5$. The first layer of the MERA coarse-grains 12 indices into 6 indices and each of the isometries in this layer is a $10 \times 10 \times 5$ tensor. The second layer of the MERA coarse-grains the remaining 6 indices of the first layer into 3 indices and therefore consists of $5 \times 5 \times 5$ isometries. The top tensor of the MERA is a 3-way cubical tensor with dimension 5. All isometries and disentanglers are initialized as random matrices, drawn from an standard normal distribution, which are then made orthogonal or orthonormal through a QR decomposition. The top tensor is also initialized as a random matrix. A comparison of the TT and MERA in terms of how well they compress the original $10^{12}$ is given in Table~\ref{tab:MERAvsTT}.
The corresponding TT has TT-ranks $R_2=10, R_3=100, R_4=50, R_5=500, R_6=250, R_7=2500, R_8=250, R_9=500, R_{10}=50, R_{11}=100, R_{12}=10$ and needs 15620200 elements. This constitutes a saving in storage space of $10^{12}/15620200= 6.40\times 10^{4}$. The MERA on the other hand consists of 54750 elements and this results in a saving of storage space of $10^{12} /54750 = 1.82\times 10^7$. The MERA is therefore about 285 times smaller as the TT.\\

Using Algorithm~\ref{alg:TT2MERA} to convert the TT back into a MERA with an identical structure as the ``true" MERA ($K=2$ and $S=5$) takes $32.74$ seconds and results in a relative approximation error of $1.00$. This large approximation error is explained by the truncated HOSVD (line 9 in Algorithm~\ref{alg:TT2MERA}) step not being able to truncate the ranks without losing accuracy. Using Algorithm~\ref{alg:TT2MERA} to convert the TT back into a MERA and using Algorithm~\ref{alg:disentangler} for the disentangler computation takes $63.81$ seconds. Setting the stopping criterion for Algorithm~\ref{alg:disentangler} to $\sigma_{R'}/\sigma_{R'+1} > 10^{13}$ guarantees that a tolerance of $10^{-12}$ can be used for the truncated HOSVD, thus obtaining a $K=2, S=5$ MERA with a relative approximation error of $1.16 \times 10^{-13}$. The low-rank approximation used in Algorithm~\ref{alg:disentangler} contained $5,25,25,25,5$ terms for the five disentanglers in the first layer, respectively, and $5$ terms for the three disentanglers in the second layer. The $63.81$ seconds run-time was dominated by Algorithm~\ref{alg:disentangler} reducing $R_7=2500$ down to a rank of 25, which took 53.88 seconds. The remaining 10 seconds were spent in the reduction of the ranks $R_6=R_8=250$ whereas the computation of all remaining tensors in the MERA took fractions of seconds.

\begin{table}[t]
\centering
\caption{Comparison of storage requirement and compression capability between a TT and a MERA for a 12th-order cubical tensor.}%\vspace{3pt}
\label{tab:MERAvsTT} % is used to refer this table in the text
\begin{tabular}{@{}lrr@{}}
                    & Storage requirement    & Compression \\\midrule
original tensor & $10^{12}$     & 1     \\
TT          & 15,620,200     & $6.40\times 10^{4}$   \\
MERA        & 54,750     & $ 1.82\times 10^7$
\end{tabular}
\end{table}

\section{Conclusions}
This article has introduced two new algorithms for the conversion of a TT into a Tucker decomposition and a MERA. The computation of a MERA-layer was shown to consist of one HOSVD-step for the computation of the disentanglers and one truncated HOSVD-step for the computation of the isometries. Using HOSVD to compute disentanglers was shown to be sub-optimal in terms of reducing the rank and an iterative orthogonal Procrustes algorithm was proposed that is able to find rank-lowering disentanglers. Numerical experiments have demonstrated the efficacy of the proposed algorithms. The TT to Tucker decomposition algorithm was demonstrated to be fast compared to the conventional HOSVD algorithm and resulted in an improvement of storage complexity that was one order of magnitude smaller. The MERA was shown to have even more potential in storage complexity in an experiment involving a tensor that consisted of $10^{12}$ elements where a compression improvement of a factor 285 compared to a TT was observed. The effectiveness and limitations of the orthogonal Procrustes algorithm were also explored in numerical experiments. The exact conditions under which this orthogonal Procrustes converges to a disentangler that retrieves an exact minimal-rank solution is still a topic for future research.

%\section*{Acknowledgments}
%We would like to acknowledge the assistance of volunteers in putting together this example manuscript and supplement.

\section*{Conflict of interest statement}
On behalf of all authors, the corresponding author states that there is no conflict of interest.

\bibliographystyle{siamplain}
\bibliography{references}
\end{document}